\documentclass[journal]{IEEEtran}

\usepackage{cite}
\usepackage{amsmath,amssymb,amsfonts,bm}
\usepackage{bbm}
\usepackage{graphicx}
\usepackage{subfigure}
\usepackage{textcomp}
\usepackage{algorithm}
\usepackage{algpseudocode}

\makeatletter
\renewcommand*\env@matrix[1][*\c@MaxMatrixCols c]{%
	\hskip -\arraycolsep
	\let\@ifnextchar\new@ifnextchar
	\array{#1}}
\makeatother
\def \qed {\hfill \vrule height6pt width 6pt depth 0pt}
\usepackage{mathrsfs}
\usepackage{xcolor}
\allowdisplaybreaks[4]
\newtheorem{definition}{Definition}
\newtheorem{theorem}{Theorem}
\newtheorem{lemma}{Lemma}
\newtheorem{assum}{Assumption}

\newtheorem{remark}{Remark}

\newenvironment{proof}{\begin{IEEEproof}}{\end{IEEEproof}}

\begin{document}
	\title{Data-Enabled Policy Optimization for Direct Adaptive Learning of the LQR
		\thanks{Research of F. Zhao and K. You was supported by National Science and Technology Major Project of China (2022ZD0116700) and National Natural Science Foundation of China (62033006, 62325305). (Corresponding author: Keyou You)}
		\thanks{F. Zhao and K. You are with the Department of Automation and BNRist, Tsinghua University, Beijing 100084, China. (e-mail: zhaofr18@tsinghua.org.cn, youky@tsinghua.edu.cn)}\thanks{F. D\"{o}rfler is with the Department of Information Technology and Electrical Engineering, ETH Z\"{u}rich, 8092 Z\"{u}rich, Switzerland. (e-mail: dorfler@control.ee.ethz.ch)}\thanks{A. Chiuso is with the Department of Information Engineering, University of Padova, Via Gradenigo 6/b, 35131 Padova, Italy. (e-mail: alessandro.chiuso@unipd.it)}
		\author{Feiran Zhao, Florian D\"{o}rfler, Alessandro Chiuso, Keyou You}}
	\maketitle
	
	\begin{abstract}
		Direct data-driven design methods for the linear quadratic regulator (LQR) mainly use offline or episodic data batches, and their online adaptation remains unclear. In this paper, we propose a direct adaptive method to learn the LQR from online closed-loop data. First, we propose a new policy parameterization based on the sample covariance to formulate a direct data-driven LQR problem, which is shown to be equivalent to the certainty-equivalence LQR with optimal non-asymptotic guarantees. Second, we design a novel data-enabled policy optimization (DeePO) method to directly update the policy, where the gradient is explicitly computed using only a batch of persistently exciting (PE) data. Third, we establish its global convergence via a projected gradient dominance property. Importantly, we efficiently use DeePO to adaptively learn the LQR by performing only one-step projected gradient descent per sample of the closed-loop system, which also leads to an explicit recursive update of the policy. Under PE inputs and for bounded noise, we show that the average regret of the LQR cost is upper-bounded by two terms signifying a sublinear decrease in time $\mathcal{O}(1/\sqrt{T})$ plus a bias scaling inversely with signal-to-noise ratio (SNR), which are independent of the noise statistics. Finally, we perform simulations to validate the theoretical results and demonstrate the computational and sample efficiency of our method.
	\end{abstract}
	
	\begin{IEEEkeywords}
		Adaptive control, linear quadratic regulator, policy optimization, direct data-driven control.
	\end{IEEEkeywords}
	
	\section{Introduction}
	
\begin{figure}[t]
		\centerline{\includegraphics[height=17mm]{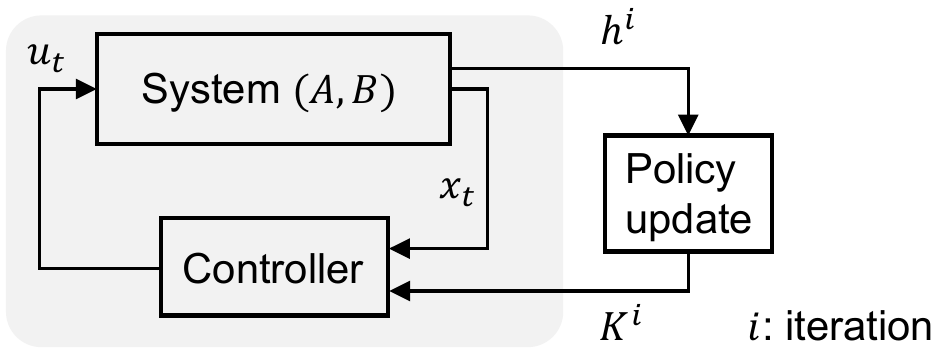}}
		\caption{An illustration of episodic approaches, where $h^i=(x_0,u_0,\dots,x_{T^i})$ denotes the $i$-th episode of data, {and the episodes can be consecutive.} }
		\label{subfig:episodic}
	\end{figure}

\begin{figure}[t]
		\centerline{\includegraphics[height=25mm]{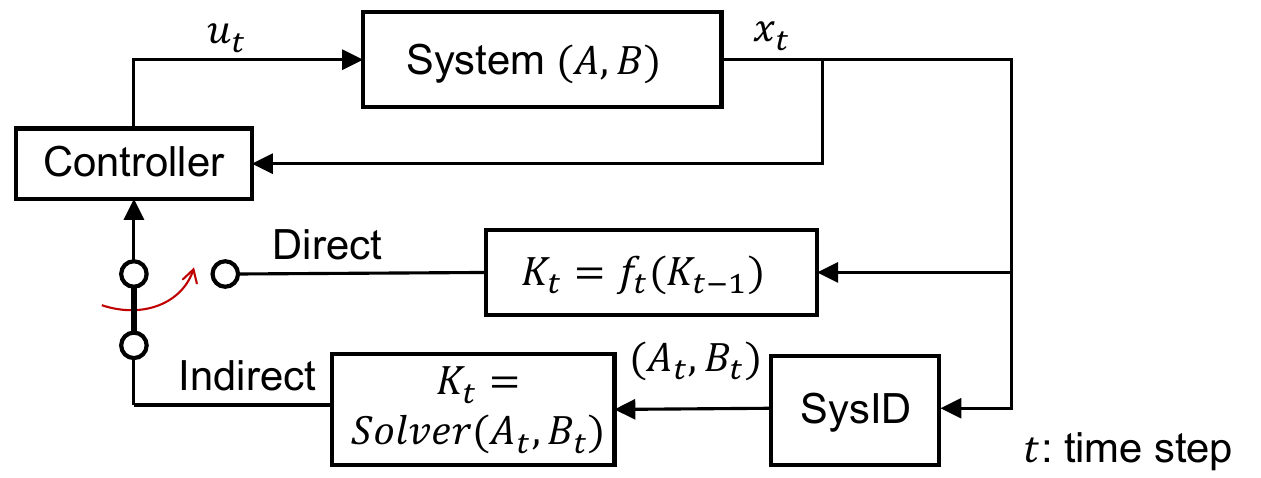}}
		\caption{An illustration of indirect and direct adaptive approaches in a closed-loop system, where function $f_t$ has an explicit form.}
		\label{subfig:adaptive}
	\end{figure}

	
	As a cornerstone of modern control theory, the linear quadratic regulator (LQR) design has been widely studied in data-driven control, where no model but only raw data is available~\cite{anderson2007optimal}. The manifold approaches to data-driven LQR design can be broadly categorized as \textit{indirect}, i.e., based on system identification (SysID) followed by model-based control design, versus \textit{direct} when bypassing the identification step. {Another classification is \textit{episodic} when obtaining the control policy from one episode of data or by alternating episodes of data collection and control (see Fig. \ref{subfig:episodic}),} versus \textit{adaptive} when updating the control policy from online closed-loop data (see Fig. \ref{subfig:adaptive}).
	  
	
	The indirect data-driven LQR design has a rich history with well-understood tools for identification and control. Representative approaches include optimism-in-face-of-uncertainty \cite{abbasi2011regret},  robust control \cite{dean2020sample}, certainty-equivalence control \cite{mania_certainty_2019, simchowitz2020naive, campi1998adaptive, wang2021exact, lu2023almost}, and adaptive dynamic programming \cite{bian2016value,lewis2009reinforcement}. Most of them are episodic in that they either estimate the system dynamics from a single episode of offline data, or update their estimate only after an episode is completed~\cite{abbasi2011regret, dean2020sample, mania_certainty_2019, simchowitz2020naive}. This is due to their requirement of statistically independent data and regret analysis methods. {Notable adaptive methods are rooted on certainty-equivalence LQR \cite{campi1998adaptive, wang2021exact, lu2023almost}}: a system is first identified by ordinary least-squares from closed-loop data, and then a certainty-equivalence LQR is obtained with Riccati equations by treating the estimated system as the ground-truth~\cite{mania_certainty_2019}. By alternating identification and certainty-equivalence LQR, they guarantee convergence to the optimal LQR gain. In particular, the work \cite{ campi1998adaptive} takes the first step towards indirect adaptive control with asymptotic convergence guarantees by regularizing the identification objective with the LQR cost. Recent works \cite{wang2021exact, lu2023almost} have shown that certainty-equivalence control with explorative input ensuring persistency of excitation meets optimal non-asymptotic guarantees.

	Different from the indirect design, direct methods entirely by-pass models; see \cite{dorfler2023data} for a discussion of the trade-offs. An emerging line of direct methods obtains the LQR directly from a single episode of persistently exciting (PE) data~\cite{celi2023closed,de2019formulas,de2021low,dorfler2021certainty,dorfler22on,10549788}. It is inspired by subspace methods~\cite{markovsky2023data} and the fundamental lemma~\cite{willems2005note} in behavioral system theory \cite{markovsky2021behavioral, coulson2019data, breschi_data-driven_2023,chiuso2023harnessing}. {Using subspace relations, the works \cite{celi2023closed,de2019formulas,de2021low} show that the closed-loop system can be parameterized by state-space data, leading to direct data-driven formulations of the LQR problem.} By a change of variables, they can be reformulated as semi-definite programs (SDPs) parameterized by raw data matrices. In the presence of noise, regularization is introduced for direct LQR design to promote certainty-equivalence or robustness~\cite{de2021low,dorfler2021certainty,dorfler22on}. There are also works~\cite{van2020noisy, 9903316} developing matrix S-lemma or combining prior knowledge for robust LQR design, while they are inherently conservative. {Though these methods only use a single episode of offline data, the dimension of their formulations usually scales with the data length. Since adaptation of their policy to the latest data may not improve control performance, they cannot use online closed-loop data to achieve adaptive learning of the LQR.} In fact, their real-time adaptation is acknowledged as an open problem in the data-driven control field~\cite{markovsky2021behavioral,markovsky2023data,annaswamy2023adaptive}.
	
	A potential path towards direct and online adaptive control is policy optimization (PO), a direct design framework where the policy is parameterized and recursively updated to minimize a cost function. Dating back to the adaptive control of aircraft in the 1950s \cite{whitaker1958design}, the concept of direct PO has a long history in control theory~\cite{1104686, kalman1960contributions, maartensson2009gradient}. However, due to the non-convexity of PO formulations, it is usually challenging to obtain strong performance guarantees. Recently, there have been resurgent interests in studying theoretical properties of zeroth-order PO, which is also an essential approach of modern reinforcement learning~\cite{fazel2018global,mohammadi2022convergence, malik2019derivative,zhao23global,zhao2022sample,bin2022towards}. It improves the policy by gradient methods, where the gradient is estimated from observations of the cost. For the LQR learning problem, zeroth-order PO meets global linear convergence thanks to a gradient dominance property~\cite{fazel2018global,mohammadi2022convergence, malik2019derivative}. However, zeroth-order PO is intrinsically unsuitable for adaptive control since (a) the cost used for gradient estimate can be obtained only \textit{after} observing an entire trajectory, (b) the trajectory needs to be sufficiently long to reduce the estimation error, and (c) it requires numerous trajectories to find an optimal policy. Different from zeroth-order PO, our recent work~\cite{zhao2023data} has proposed a data-enabled PO (DeePO) method for the LQR, where the gradient is computed directly from a single trajectory of finite length, and shown global convergence. This is achieved by adopting the data-based policy parameterization in \cite{de2019formulas,de2021low,dorfler2021certainty,dorfler22on}. While this DeePO method is based on offline data, it paves the way to applying PO for direct and online adaptive control.


	Following \cite{zhao2023data}, this paper proposes a novel DeePO method for direct adaptive learning of the LQR, where the policy is directly updated by gradient methods from online closed-loop data; see Fig. \ref{subfig:adaptive} for an illustration. Hence, we provide a promising solution to the open problem in~\cite{markovsky2021behavioral,markovsky2023data,annaswamy2023adaptive}. Our contributions are summarized below.
	\begin{itemize} 
	 \item We propose a new policy parameterization for the LQR based on sample covariance, which is a key ingredient in SysID~\cite{mania_certainty_2019}, filtering~\cite{kailath2000linear}, and data-driven control parameterizations~\cite{rantzer2024data, rantzer2021minimax, chiuso2023harnessing,song2024role}. Compared with the existing parameterization \cite{de2019formulas,de2021low,dorfler2021certainty,dorfler22on}, it has two salient features that enable online adaptation of DeePO. First, the dimension of the parameterized policy remains constant depending only on the system dimension. Second, the resulting direct LQR is shown to be equivalent to the indirect certainty-equivalence LQR, which usually requires regularized formulations and methods \cite{de2021low,dorfler2021certainty,dorfler22on}. In view of \cite{mania_certainty_2019,simchowitz2020naive}, this equivalence implies that our covariance parameterization enables sample-efficient online learning. The covariance parameterization can also be used to solve other control problems (e.g., the $\mathcal{H}_{\infty}$ problem \cite{de2019formulas}) in a direct data-driven fashion. 
	  
	  \item We propose a DeePO method to solve the covariance-parameterized LQR problem with offline data and show global convergence. The key to our analysis is a \textit{projected gradient} dominance property, which is distinguished from the usual \textit{gradient} dominance in PO literature \cite{fazel2018global,mohammadi2022convergence, malik2019derivative}. 
	  
	  \item {We use DeePO to adaptively learn the optimal LQR gain from online closed-loop data, starting from an initial stabilizing gain learned from offline data.} The proposed approach is direct, online, and has an explicit recursive update of the policy. Moreover, it can be extended straightforwardly to time-varying systems by adding a forgetting factor to the covariance parameterization.
 
 	\item We provide non-asymptotic guarantees of DeePO for adaptive learning of the LQR, which are independent of noise statistics. {Our focus is on the convergence of the policy instead of stability, which is in line with the  RL perspective \cite{lewis2009reinforcement,fazel2018global} and the definition of adaptive control by Zames \cite{zames1998adaptive}.} Under  PE inputs and bounded noise, we show that the average regret of the LQR cost is upper-bounded by two terms signifying a sublinear decrease in time $\mathcal{O}(1/\sqrt{T})$ plus a bias scaling inversely as signal-to-noise ratio (SNR). This convergence result improves over single batch methods~\cite{de2021low,dorfler2021certainty,dorfler22on}, whose performance also depends on SNR but does not decay over time. Our sublinear decrease rate aligns with that of first-order methods in online convex optimization of smooth functions~\cite{hazan2016introduction}, even though our considered LQR problem is non-convex. This sublinear rate shows the sample efficiency of DeePO to learn from online data.   
	 \end{itemize} 
	

	In the simulations, we validate the global convergence of DeePO. Moreover, we compare DeePO with the indirect adaptive approach~\cite{ campi1998adaptive, wang2021exact, lu2023almost} and zeroth-order PO \cite{fazel2018global,malik2019derivative,mohammadi2022convergence} for  the benchmark problem in \cite{dean2020sample}. The simulations demonstrate favorable computational and sample efficiency of DeePO.

	The rest of this paper is organized as follows. Section \ref{sec:prob} recapitulates data-driven LQR formulations and describes the adaptive learning problem. Section \ref{sec:newpara} proposes the covariance parameterization for the LQR. Section \ref{sec:deepo} uses DeePO to solve LQR parameterization with offline data. Section \ref{sec:adap} extends DeePO to direct adaptive learning of the LQR with non-asymptotic guarantees. Section \ref{sec:simu} validates our theoretical results via simulations. We provide concluding remarks in Section \ref{sec:conclution}. All the proofs are deferred to   Appendices.

	\textbf{Notation.} We use $I_n$ to denote the $n$-by-$n$ identity matrix. We use $\underline{\sigma}(\cdot)$ to denote the minimal singular value of a matrix. We use $\|\cdot\|$ to denote the $2$-norm of a vector or matrix, and $\|\cdot\|_F$ the Frobenius norm. We use $\rho(\cdot)$ to denote the spectral radius of a square matrix. We use $A^\dagger:=A^{\top}(AA^{\top})^{-1}$ to denote the right inverse of a full row rank matrix $A\in \mathbb{R}^{n\times m}$. We use $\mathcal{N}(A)$ to denote the nullspace of $A$ and $\Pi_A := I_{m}-A^{\dagger}A$ the projection operator onto $\mathcal{N}(A)$.


	\section{Data-driven formulations and adaptive learning of the LQR}\label{sec:prob}
	In this section, we first recapitulate indirect certainty-equivalence LQR (least-square SysID followed by model-based LQR design) \cite{dorfler22on}, direct LQR design using data-based policy parameterization \cite{de2019formulas,de2021low}, and policy optimization (PO) of the LQR based on zeroth-order gradient estimate~\cite{fazel2018global}. Then, we formalize our adaptive learning problem.
	\subsection{Model-based LQR}

	Consider a discrete-time LTI system
	\begin{equation}\label{equ:sys}
	\left\{\begin{aligned}
	x_{t+1} & =A x_t+B u_t+w_t \\
	z_t & =\begin{bmatrix}
	Q^{1 / 2} & 0 \\
	0 & R^{1 / 2}
	\end{bmatrix}
	\begin{bmatrix}
	x_t \\
	u_t
	\end{bmatrix}
	\end{aligned}\right.,
	\end{equation}
	where $t\in \mathbb{N}$, $x_t\in\mathbb{R}^{n}$ is the state, $u_t\in\mathbb{R}^{m}$ is the control input, $w_t \in \mathbb{R}^n$ is the noise, and $z_t$ is the performance signal of interest. We assume that $(A,B)$ are controllable and the weighting matrices $(Q, R)$ are positive definite. {Throughout the paper, we do not assume any  statistics of the noise $w_t$.}
	
	The LQR problem is phrased as finding a state-feedback gain $K\in \mathbb{R}^{m\times n}$ that minimizes the $\mathcal{H}_2$-norm of the transfer function $\mathscr{T}(K):w \rightarrow z$ of the closed-loop system
	\begin{equation}\label{equ:closedsys}
	\begin{bmatrix}
	x_{t+1} \\
	z_t
	\end{bmatrix}=\begin{bmatrix}[c|c]
	A+BK & I_n \\
	\hline \begin{bmatrix}
	Q^{1 / 2} \\
	R^{1 / 2} K
	\end{bmatrix} & 0
	\end{bmatrix}\begin{bmatrix}
	x_t \\
	w_t
	\end{bmatrix}.
	\end{equation}
	When $A+BK$ is stable, it holds that \cite{anderson2007optimal}
	\begin{equation}\label{equ:transfer}
	\|\mathscr{T}(K)\|_2^2  = \text{Tr}((Q+K^{\top}RK)\Sigma_K)=:C(K),
	\end{equation}
	where $\Sigma_K$ is the closed-loop state covariance matrix obtained as the positive definite solution to the Lyapunov equation
	\begin{equation}\label{equ:Sigma}
	\Sigma_K = I_n + (A+BK)\Sigma_K (A+BK)^{\top}.
	\end{equation}
	We refer to $C(K)$ as the LQR cost and to (\ref{equ:transfer})-(\ref{equ:Sigma}) as a \textit{policy parameterization} of the LQR.
	
	There are alternative formulations to find the optimal LQR gain $K^*:=\arg\min_{K}C(K)$ of (\ref{equ:transfer})-(\ref{equ:Sigma}), e.g., via the celebrated Riccati equation  with known $(A,B)$ \cite{anderson2007optimal}. For unknown $(A,B)$, there is a plethora of data-driven control methods to find $K^*$, some of which we recapitulate below.
	
	
	%

	\subsection{Indirect certainty-equivalence LQR with ordinary least-square identification}
	This conventional approach to data-driven LQR design follows the certainty-equivalence principle: it first identifies a system $(A,B)$ from data, and then solves the LQR problem regarding the identified model as the ground-truth. The SysID step is based on subspace relations among the state-space data. Consider a $t$-long time series\footnote{The time series do not have to be consecutive. All results in Sections \ref{sec:prob}-\ref{sec:deepo} also hold when each column in (\ref{equ:data}) is obtained from independent experiments or even averaged data sets.} of states, inputs, noises, and successor states
	\begin{equation}\label{equ:data}
	\begin{aligned}
	X_{0} &:= \begin{bmatrix}
	x_0& x_1& \dots& x_{t-1}
	\end{bmatrix}\in \mathbb{R}^{n\times t},\\
	U_{0} &:= \begin{bmatrix}
	u_0& u_1& \dots& u_{t-1}
	\end{bmatrix}\in \mathbb{R}^{m\times t}, \\
	W_{0} &:= \begin{bmatrix}
	w_0& w_1& \dots& w_{t-1}
	\end{bmatrix}\in \mathbb{R}^{n\times t}, \\
	X_{1} &:= \begin{bmatrix}
	x_1& x_2& \dots& x_t
	\end{bmatrix}\in \mathbb{R}^{n\times t},
	\end{aligned}
	\end{equation}
	which satisfy the system dynamics
	\begin{equation}\label{equ:dynamics}
	X_1 = AX_0+ BU_0 + W_0.
	\end{equation}
	
	We assume that the data is {\em persistently exciting (PE)} \cite{willems2005note}, i.e., the block matrix of input and state data 
	$$
	D_0 := \begin{bmatrix}
	U_0 \\
	X_0
	\end{bmatrix}\in \mathbb{R}^{(m+n)\times t}
	$$
	has full row rank
	\begin{equation}\label{equ:rank}
	\text{rank}(D_0) = m+n.
	\end{equation}
	This PE condition is necessary for the data-driven LQR design~\cite{van2020data,kang2023minimum}.
	
	Based on the subspace relations \eqref{equ:dynamics} and the rank condition (\ref{equ:rank}), an estimate $(\widehat{A},\widehat{B})$ of the system can be obtained as the unique solution to the ordinary least-squares problem
	\begin{equation}\label{equ:sysid}
	[\widehat{B}, \widehat{A}] = \underset{B, A}{\arg \min }\left\|X_1-[B,A] D_0\right\|_F = X_1D_0^{\dagger}.
	\end{equation}
	Following the certainty-equivalence principle~\cite{dorfler2021certainty}, the system $(A,B)$ is replaced with its estimate $(\widehat{A},\widehat{B})$ in (\ref{equ:transfer})-(\ref{equ:Sigma}), and the LQR problem can be reformulated as a bi-level program
	\begin{equation}\label{prob:indirect}
	\begin{aligned}
	\mathop{\text{minimize}}\limits_{K, \Sigma\succeq 0} ~~&\text{Tr}\left((Q+K^{\top}RK)\Sigma\right),\\
	\text{subject to} ~~&\Sigma = I_n + (\widehat{A} + \widehat{B}K)\Sigma (\widehat{A} + \widehat{B}K)^{\top}\\
	&[\widehat{B}, \widehat{A}] = \underset{B, A}{\arg \min }\left\|X_1-[B,A] D_0\right\|_F.
	\end{aligned}
	\end{equation}
	The problem \eqref{prob:indirect} is termed {\em certainty-equivalence} and {\em indirect data-driven} LQR design in \cite{dorfler22on}. By alternating between solving \eqref{prob:indirect} and using the resulting policy to collect online closed-loop data, it achieves adaptive learning of the LQR with optimal non-asymptotic guarantees~\cite{wang2021exact, lu2023almost}. 
	

	\subsection{Direct LQR with data-based policy parameterization}
	In contrast to the SysID-followed-by-control approach \eqref{prob:indirect}, the{ \em direct data-driven} LQR design aims to find $K^*$ bypassing the identification step \eqref{equ:sysid}~\cite{de2019formulas,de2021low,dorfler2021certainty}, which we recapitulate as follows. It uses a data-based policy parameterization: by the rank condition \eqref{equ:rank}, there exists a matrix $G \in \mathbb{R}^{T\times n}$ that satisfies
	\begin{equation}\label{equ:relation}
	\begin{bmatrix}
	K \\
	I_n
	\end{bmatrix}=
	D_0G
	\end{equation}
	for any given $K$. Together with the subspace relation (\ref{equ:dynamics}), the closed-loop matrix can be written as
	$$
	A+BK=[B,A]\begin{bmatrix}
	K \\
	I_n
	\end{bmatrix}\overset{\eqref{equ:relation}}{=}[B,A]D_0G\overset{\eqref{equ:dynamics}}{=}(X_1 - W_0)G,
	$$
	Since $W_0$ is unknown and unmeasurable, it is disregarded and $X_1G$ is used as the closed-loop matrix. Following the certainty-equivalence principle~\cite{dorfler2021certainty}, we substitute $A+BK$ with $X_1G$ in (\ref{equ:transfer})-(\ref{equ:Sigma}), and together with \eqref{equ:relation} the LQR problem becomes 
	\begin{equation}\label{prob:equi}
	\begin{aligned}
	&\mathop{\text {minimize}}\limits_{G, \Sigma\succeq 0}~ \text{Tr}\left((Q+G^{\top}U_0^{\top}RU_0G)\Sigma\right),\\
	&\text{subject to}~~ \Sigma = I_n + X_1G\Sigma G^{\top}X_1^{\top},~X_0G =I_n 
	\end{aligned}
	\end{equation}
	with the gain matrix $K = U_0G$, which can be reformulated as an SDP~\cite{de2019formulas}. The LQR parameterization \eqref{prob:equi} is direct data-driven, as it does not involve any explicit SysID. 
	
	{However, the LQR parameterization (\ref{prob:equi}) is not robust to noise and requires regularization~\cite{dorfler2021certainty}.} Moreover, its dimension scales linearly with $t$, and it is unclear how to turn (\ref{prob:equi}) into a recursive formulation for online closed-loop adaptation. {These issues are naturally addressed with the covariance parameterization in the sequel.}


	\subsection{PO of the LQR using zeroth-order gradient estimate}
	As an essential approach of modern reinforcement learning, zeroth-order PO \cite{fazel2018global,malik2019derivative,mohammadi2022convergence} finds $K^*$ by iterating the gradient descent  for the policy parameterization of the LQR (\ref{equ:transfer})-(\ref{equ:Sigma}):
	\begin{equation}\label{equ:pogd}
	K' = K - \eta \widehat{\nabla C(K)}.
	\end{equation}
	The PO \eqref{equ:pogd} is initialized with a stabilizing gain,  $\eta>0$ is a constant stepsize, and $\widehat{\nabla C(K)}$ is the gradient estimate from zeroth-order information, e.g., the two-point gradient estimate 
	\begin{equation}\label{equ:gdest}
	\widehat{\nabla C(K)} = \left(\widehat{C}(K+ rF)-\widehat{C}(K-rF)\right) \cdot \frac{mn}{r} F,
	\end{equation}
	where $r>0$ is the smoothing radius, $F$ is uniformly sampled from the unit sphere $\mathbb{S}^{mn-1}$, and {$\widehat{C}$ is the approximated cost observed from a single realization of a $T$-long trajectory of system (\ref{equ:sys}), i.e.,
	$
	\widehat{C}(K)= \frac{1}{T}\sum_{t=0}^{T-1}\|z_t\|^2
	$ with $u_t = Kx_t$.}
	
	While zeroth-order PO has a recursive policy update, it requires to sample a large number of long trajectories to find $K^*$, which is inherently unsuitable for online adaptive control.
	
	\subsection{Direct adaptive learning for the LQR with online closed-loop data}
	All the aforementioned data-driven LQR design methods either lack a recursive policy update, or are unsuitable for adaptive control with online closed-loop data. In this paper, we focus on the  direct adaptive learning problem of the LQR:
	
	\textbf{Problem:} design a recursive direct method based on online closed-loop data such that the control policy converges to the optimal LQR gain.
	
	{Our objective concerns the optimality of the policy whereas traditional adaptive control focuses on stability. Our perspective is in line with the RL perspective \cite{lewis2009reinforcement,fazel2018global} and the essence of adaptive control by Zames \cite{zames1998adaptive}, namely, improving over the best control with \textit{a prior} information.} Following \cite{de2019formulas,de2021low,dorfler2021certainty,dorfler22on}, we start with a direct data-driven LQR formulation but with a new policy parameterization. Then, we take an iterative PO perspective to the direct LQR as in our previous work \cite{zhao2023data}, and use gradient methods to achieve adaptive learning of the LQR with online closed-loop data.

	\section{Direct data-driven LQR with a new policy parameterization}\label{sec:newpara}
	In this section, we first propose a new policy parameterization based on the sample covariance  to formulate the direct data-driven LQR problem. Then, we establish its equivalence to the indirect certainty-equivalence LQR problem.
	
	\subsection{A new policy parameterization using sample covariance}
	To efficiently use data, we propose a new policy parameterization based on the sample covariance of input-state data
	\begin{equation}\label{def:cov}
	\Phi :  = \frac{1}{t}D_0D_0^{\top}=
	\begin{bmatrix}
	U_0D_0^{\top}/t  \\
	 \hline X_0D_0^{\top}/t 
	\end{bmatrix}
	=\begin{bmatrix}
	\overline{U}_0 \\
 \hline	\overline{X}_0
	\end{bmatrix},
	\end{equation}
	which plays an important role in SysID~\cite{mania_certainty_2019}, filtering~\cite{kailath2000linear}, and data-driven control parameterizations~\cite{rantzer2021minimax, rantzer2024data, chiuso2023harnessing,song2024role}. Under the PE rank condition (\ref{equ:rank}), the sample covariance $\Phi$ is positive definite, and there exists a \textit{unique} solution $V\in \mathbb{R}^{(n+m)\times n}$ to
	\begin{equation}\label{equ:newpara}
	\begin{bmatrix}
	K \\
	I_n
	\end{bmatrix}= \Phi V
	\end{equation}
	for any given $K$. We refer to \eqref{equ:newpara} as the \textit{covariance parameterization} of the policy. In contrast to the  parameterization in (\ref{equ:relation}), the dimension of $V$ is independent of the data length. 
	
	\begin{remark}
		The covariance parameterization (\ref{equ:newpara}) can also be derived from  (\ref{equ:relation}).  Since $D_0$ does not necessarily have full column rank, there is a considerable nullspace in the solution of (\ref{equ:relation}), which can be reparameterized via orthogonal decomposition:
		\begin{equation}\label{equ:rela}
		G = \frac{1}{t} D_0^{\top}V + \Delta, \Delta \in \mathcal{N}(D_0) 
		\end{equation}
		with $V \in \mathbb{R}^{(n+m)\times n}$. If we remove the nullspace $\mathcal{N}(D_0)$ in (\ref{equ:rela}), then the parameterization (\ref{equ:relation}) reduces to (\ref{equ:newpara}). Note that the nullspace is undesirable, and regularization methods are usually used  to single out a favorable solution\cite{de2021low,dorfler2021certainty,dorfler22on}. \qed
	\end{remark}

	For brevity, define the sample covariances
	$
	\overline{W}_0= W_0D_0^{\top}/t, \overline{X}_1= X_1D_0^{\top}/t,
	$ {whose dimensions do not depend on $t$.}
	Then, the closed-loop matrix can be written as
	\begin{equation}\label{equ:clos}
	A+BK=[B,A]\begin{bmatrix}
	K \\
	I_n
	\end{bmatrix}\overset{\eqref{equ:newpara}}{=}[B,A]\Phi V\overset{\eqref{equ:dynamics}}{=}(\overline{X}_1 - \overline{W}_0)V.
	\end{equation}
	Analogous to \eqref{prob:equi}, we disregard the uncertainty  $\overline{W}_0$ in the parameterized closed-loop matrix and formulate the direct data-driven LQR problem using $(X_0,U_0,X_1)$ as
	\begin{equation}\label{prob:equiV}
	\begin{aligned}
	&\mathop{\text {minimize}}\limits_{V, \Sigma_V\succeq 0}~J(V) :=\text{Tr}\left((Q+V^{\top}\overline{U}_0^{\top}R\overline{U}_0V)\Sigma_V\right),\\
	&\text{subject to}~ ~\Sigma_V = I_n + \overline{X}_1V\Sigma_V V^{\top}\overline{X}_1^{\top},\overline{X}_0V= I_n
	\end{aligned}
	\end{equation}
	with the gain matrix $K = \overline{U}_0V$. We refer to (\ref{prob:equiV}) as the LQR problem with covariance parameterization. 
	{\begin{remark}
		In comparison with \eqref{prob:equi}, the covariance parameterization can well display and partially mitigate the effects of stochastic noise. Let $\{w_t\}$ be additive white noise. Then, it follows that $\mathbb{E}[\sum_{i=0}^{t-1}w_iu_i^{\top}]/t=0$ and $\mathbb{E}[\sum_{i=0}^{t-1}w_ix_i^{\top}]/t=0$. Approximating the expectation using data then leads to $\overline{W}_0 \approx 0$ and hence $A+BK \approx \overline{X}_1V$. A detailed analysis of the stochastic setting is outside the scope of this paper and left to future work. \qed
	\end{remark}}
	
	
	\subsection{The equivalence between the covariance parameterization of the LQR and the indirect certainty-equivalence LQR}
	We now show that the data-driven covariance parameterization \eqref{prob:equiV} is equivalent to the indirect certainty-equivalence LQR (\ref{prob:indirect}) in the sense that their solutions coincide. Let $J^*$ and $C_{\text{CE}}^*$ be the optimum of \eqref{prob:equiV} and (\ref{prob:indirect}), respectively. Then, we have the following result.
	\begin{lemma}[\textbf{Equivalence to the certainty-equivalence LQR}]\label{lem:equi}
		The feasible sets of \eqref{prob:indirect} and \eqref{prob:equiV} coincide under the change of variables $V = \Phi^{-1}[K^{\top},I_n]^{\top}$, and $J^* = C_{\text{CE}}^*$. 
	\end{lemma}
	
	\begin{remark}[\textbf{Implicit regularization}]\label{remark}
		This equivalence has also been shown for the LQR parameterization~\eqref{prob:equi} with a sufficiently large certainty-equivalence regularizer, i.e.,
		\begin{equation}\label{equ:Gce}
		\begin{aligned}
		&\mathop{\text {minimize}}\limits_{G, \Sigma\succeq 0}~ \text{Tr}\left((Q+G^{\top}U_0^{\top}RU_0G)\Sigma\right) + \lambda \|\Pi_{D_0}G \|,\\
		&\text{subject to}~~ \Sigma = I_n + X_1G\Sigma G^{\top}X_1^{\top},~X_0G =I_n
		\end{aligned}
		\end{equation}
		with the gain $K = U_0G$~\cite[Corollary 3.2]{dorfler2021certainty}. Thus,  the covariance parametrization of the LQR (\ref{prob:equiV}) is \textit{implicitly} regularized as it achieves this equivalence without any regularization. An implicit regularization property is also established in the PO of \eqref{equ:Gce} in the absence of noise \cite[Theorem 2]{zhao2023data}. \qed
	\end{remark}

	While the LQR problem with covariance parameterization (\ref{prob:equiV}) can be reformulated into an SDP as in \cite{de2019formulas,de2021low,dorfler2021certainty,dorfler22on}, our subsequent solution will not make use of that. Instead, the next section solves (\ref{prob:equiV}) via an iterative PO method with convergence guarantees.

	
	\section{DeePO for the LQR with covariance parameterization using offline data}\label{sec:deepo}
	In this section, we first present our novel PO method to solve the LQR problem with covariance parameterization (\ref{prob:equiV}) given offline data $(X_0,U_0,X_1)$. Then, we show its global convergence by proving a projected gradient dominance property.
	\subsection{Data-enabled policy optimization to solve the LQR problem with covariance parameterization (\ref{prob:equiV})}
	
	We assume the feasibility of (\ref{prob:equiV}).
	\begin{assum}\label{assumption}
		The covariance parameterization of the LQR problem in (\ref{prob:equiV}) has at least one solution, i.e., the feasible set  
		$\mathcal{S}:=\{V\mid \overline{X}_0V =I_n,  \rho (\overline{X}_1V)<1\}$ is non-empty.
	\end{assum}
	\begin{remark}
		{By Lemma \ref{lem:equi}, Assumption \ref{assumption} is equivalent to the feasibility of the certainty-equivalence LQR problem \eqref{prob:indirect}, and hence is necessary also for indirect control \cite{wang2021exact}. \qed }
	\end{remark}
 
	We first derive the gradient for $J(V)$. Since $\Sigma_V = \sum_{i=0}^{\infty} (\overline{X}_1V)^i(V^{\top}\overline{X}_1^{\top})^i $, the cost $J(V)$  can be expressed as	
	\begin{align}
	J(V)&= \text{Tr}\left((Q+V^{\top}\overline{U}_0^{\top}R\overline{U}_0V)\sum_{i=0}^{\infty} (\overline{X}_1V)^i (V^{\top}\overline{X}_1^{\top})^i\right) \notag\\
	&=\text{Tr}\left(\sum_{i=0}^{\infty} (V^{\top}\overline{X}_1^{\top})^i(Q+V^{\top}\overline{U}_0^{\top}R\overline{U}_0V) (\overline{X}_1V)^i\right) \notag \\
	&=\text{Tr}(P_V), \label{def:JG}
	\end{align}
	{where $P_V \succ 0$ is the unique solution to the Lyapunov equation}
	\begin{equation}\label{equ:Lya_P}
	P_V = Q + V^{\top}\overline{U}_0^{\top}R\overline{U}_0V + V^{\top}\overline{X}_1^{\top}P_V\overline{X}_1V.
	\end{equation}
	Then, the closed-form expression for $\nabla J(V)$ is given as follows.
	\begin{lemma}\label{lem:gradient}
		For $V\in \mathcal{S}$, the gradient of $J(V)$ is
		$$
		\nabla J(V) = 2 \left(\overline{U}_0^{\top}R\overline{U}_0+\overline{X}_1^{\top}P_V\overline{X}_1\right)V \Sigma_V.
		$$
	\end{lemma}
	
	{By \eqref{prob:equiV}, \eqref{equ:Lya_P}, and Lemma \ref{lem:gradient}, the gradient can be computed directly from raw data matrices and the solution to \eqref{equ:Lya_P}.} Starting from a feasible $V^0 \in \mathcal{S}$, the projected gradient method iterates as
	\begin{equation}\label{equ:gd}
	V^{k+1} = V^k - \eta\Pi_{\overline{X}_0} \nabla J(V^k) , ~k\in\mathbb{N},
	\end{equation}
	where $\eta \geq 0$ is a constant stepsize, and the projection $\Pi_{\overline{X}_0}$ enforces the subspace constraint $\overline{X}_0V =I_n$. We refer to (\ref{equ:gd}) as \textit{\textbf{d}ata-\textbf{e}nabl\textbf{e}d \textbf{p}olicy \textbf{o}ptimization (DeePO)} as the projected gradient is data-driven, and the control gain can be recovered from the covariance parameterization (\ref{equ:newpara}) as $K = \overline{U}_0V$. 
	
	Due to non-convexity of the considered problem \eqref{def:JG}-\eqref{equ:Lya_P}, it is challenging to provide convergence guarantees for DeePO in (\ref{equ:gd}). The existing PO methods achieve global convergence for the LQR by leveraging a gradient dominance property\cite{fazel2018global,malik2019derivative,mohammadi2022convergence}. Here, we show a similar projected gradient dominance property of $J(V)$, as a key to global convergence of (\ref{equ:gd}).
	
	\subsection{Projected gradient dominance of the LQR cost}
	We first define the projected gradient dominance property.
	\begin{definition}[\textbf{Projected gradient dominance}]\label{def:pgd}
		A continuously differentiable function $g(x): \mathbb{R}^n \rightarrow \mathbb{R}$ is projected gradient dominated of degree $p\in \{1,2\}$ over a subset $\mathcal{X}\subseteq \text{dom}(g) $ if
		$$
		g(x) - g(x^*) \leq \alpha \|\Pi_{\mathcal{X}} (\nabla g(x)) \|^p, ~\forall x \in \mathcal{X},
		$$
		where $\alpha > 0$ is a constant, $x^*$ is a minimizer of $g(x)$ over $\mathcal{X}$, and $\Pi_{\mathcal{X}}$ is a projection onto a closed convex subset of  $\mathcal{X}$.
	\end{definition}
	
	In comparison to the definition of gradient dominance~\cite[Definition 1]{bin2022towards}, Definition \ref{def:pgd} focuses on the property of the \textit{projected} gradient, meaning that optimality is achieved when the projected gradient attains zero. It leads to global convergence of projected gradient descent in conjuncture with a smoothness property, and the convergence rate depends on the degree $p$. In particular, for smooth objective  $p=1$ leads to a sublinear rate and $p=2$ leads to a linear rate. 
	
	For  problem \eqref{def:JG}-\eqref{equ:Lya_P}, we show that $J(V)$ is projected gradient dominated of degree 1 over any sublevel set
	$
	\mathcal{S}(a):= \{V\in \mathcal{S} \mid J(V)\leq a\}
	$
	with $a>0$, the proof of which leverages \cite[Theorem 1]{sun2021analysis} and a convex reparameterization.
	
	{\begin{lemma}[\textbf{Projected gradient dominance of degree 1}]\label{lem:pl}
			For $V \in \mathcal{S}(a)$, it holds that
			$$
			J(V) - J^* \leq \mu(a) \|\Pi_{\overline{X}_0} \nabla J(V)\|
			$$
			with 
			$
			\mu(a) := {4a^2}/\left({(\underline{\sigma}(Q))^{{3}/{2}}(\underline{\sigma}(R))^{{1}/{2}} \underline{\sigma}(\overline{U}_0)}\right).
			$
	\end{lemma}}

	Armed with Lemma \ref{lem:pl}, we show the global convergence of the projected gradient descent (\ref{equ:gd}) in the next subsection.

	\subsection{Global convergence of DeePO}
	As in the LQR parameterization (\ref{equ:transfer})-(\ref{equ:Sigma}), the cost $J(V)$ here tends to infinity as $V$ approaches the boundary $\partial \mathcal{S}$~\cite{fazel2018global}. Thus, it is only \textit{locally} smooth over a sublevel set.

	\begin{lemma}[\textbf{Local smoothness}]\label{lem:smooth}
		For any $V,V' \in \mathcal{S}(a)$ satisfying $V+\phi(V'-V) \in \mathcal{S}(a), \forall \phi \in [0,1]$, it holds that
		$$
		J(V') \leq J(V) +  \langle \nabla J(V), V' - V \rangle + {l(a)}\|V'-V\|^2/2,
		$$
		where the smoothness constant is given by
		$$
		l(a):= 4a^2\left(\xi(a) + a -\underline{\sigma}(Q)\right) \frac{\|\overline{X}_1\|_F^2}{\underline{\sigma}^2(Q)} +  \frac{2\xi(a) a}{\underline{\sigma}(Q)}
		$$
		with $\xi(a) := \|\overline{U}_0\|^2\|R\|+\|\overline{X}_1\|^2a$. 
	\end{lemma}
	
	While Lemmas \ref{lem:pl} and \ref{lem:smooth} are local properties over a sublevel set, we can still use them by selecting an appropriate stepsize such that the policy sequence stays in this sublevel set. For brevity, let $\mu^0:= \mu(J(V^0))$ and $l^0:= l(J(V^0))$, where $V^0\in \mathcal{S}$ is an initial policy. The convergence of DeePO in (\ref{equ:gd}) is shown as follows.
	
	\begin{theorem}[\textbf{Global convergence}]\label{thm:conv}
		For $V^0\in \mathcal{S}$ and a stepsize $\eta \in (0,1/l^0]$, the projected gradient descent (\ref{equ:gd}) leads to $V^k \in \mathcal{S}(J(V^0)), \forall k \in \mathbb{N}$. Moreover, for any $\epsilon> 0$ and
		\begin{equation}\label{equ:kk}
		k \geq \frac{2(\mu^0)^2}{\epsilon(2\eta - l^0\eta^2)},
		\end{equation}
		the update (\ref{equ:gd}) has performance bound
		$
		J(V^k)-J^* \leq \epsilon.
		$
	\end{theorem}
	
	We compare our DeePO method developed thus far with  zeroth-order PO in (\ref{equ:pogd})-(\ref{equ:gdest})  \cite{fazel2018global,mohammadi2022convergence,malik2019derivative}. Since $J(V)$ is projected gradient dominated of degree 1, we have shown sublinear convergence of projected gradient descent to $J^*$ (in fact, we observe a faster linear rate in the simulation). In comparison, $C(K)$ in \eqref{equ:transfer} is gradient dominated of degree 2 \cite{fazel2018global}, and hence zeroth-order PO meets linear convergence to $C^*:=\min_{K}C(K)$. However, they rely on zeroth-order gradient estimates, which inevitably require a large number of long trajectories. In contrast, DeePO efficiently computes the gradient from a batch of PE data based on the covariance parameterization of the LQR \eqref{prob:equiV}. This remarkable feature enables DeePO to utilize online closed-loop data for adaptive control with recursive updates, as shown in the next section.
	
	\section{DeePO for direct, adaptive, and recursive learning of the LQR with online closed-loop data}\label{sec:adap}
	In this section, we first use DeePO to adaptively and recursively learn the optimal LQR gain from online closed-loop data. Then, we provide non-asymptotic convergence guarantees for DeePO.
	\subsection{Direct adaptive learning of the LQR}
	In the adaptive control setting, we collect online closed-loop data $ x_{t},u_t,x_{t+1} $ at time $t$ that constitutes the data matrices $(X_{0,t+1}, U_{0,t+1}, X_{1,t+1})$\footnote{Since the closed-loop data grows with time in the adaptive control setting, we use $X_{0,t}, U_{0,t}, W_{0,t}, X_{1,t}$ to denote the data series of length $t$ in \eqref{equ:data}. We also add a subscript $t$ to other notations to highlight the time dependence. Under the new notations, all the results in Sections \ref{sec:prob}, \ref{sec:newpara}, \ref{sec:deepo} still hold.}. The goal is to use DeePO to learn the optimal LQR gain $K^*$ directly from online closed-loop data. Our key idea is to first use  $(X_{0,t+1}, U_{0,t+1}, X_{1,t+1})$ to perform only one-step projected gradient descent for the parameterized policy at time $t$, then use the updated policy to control the system, and repeat. The details are presented in Algorithm \ref{alg:deepo}, which is both direct and adaptive.
	
	
	\begin{algorithm}[t]
		\caption{DeePO for direct adaptive learning of the LQR}
		\label{alg:deepo}
		\begin{algorithmic}[1]
			\Require Offline data $(X_{0,t_0}, U_{0,t_0}, X_{1,t_0})$, an initial stabilizing policy $K_{t_0}$, and a stepsize $\eta$.
			\For{$t=t_0,t_0+1,\dots$}
			\State Apply $u_t$ and observe $x_{t+1}$.
			\State Given $K_{t}$, solve $V_{t+1}$ via 
			\begin{equation}\label{equ:k2v}
			V_{t+1} =\Phi_{t+1}^{-1} \begin{bmatrix}
			K_{t} \\
			I_n
			\end{bmatrix}.
			\end{equation} 
			\State Perform one-step projected gradient descent
			\begin{equation}\label{equ:pgdV}
			V_{t+1}' = V_{t+1} - \eta  \Pi_{\overline{X}_{0,t+1}} \nabla J_{t+1}(V_{t+1}),
			\end{equation}
			where the gradient $\nabla J_{t+1}(V_{t+1})$ is given by Lemma \ref{lem:gradient}.
			\State Update the control gain by 
			\begin{equation}\label{equ:KV'}
			K_{t+1} = \overline{U}_{0,t+1}V_{t+1}'.
			\end{equation}
			

			\EndFor	
		\end{algorithmic}
	\end{algorithm}
	
	{An initial stabilizing gain is obtained from offline data $(X_{0,t_0}, U_{0,t_0}, X_{1,t_0})$, i.e., $K_{t_0} = \overline{U}_{0,t_0}V_{t_0}'$ with $V_{t_0}'$  the solution to (\ref{prob:equiV}). Alternatively, $K_{t_0}$ can be found via \eqref{prob:indirect} or \eqref{equ:Gce}.}

	At the online stage $t \geq t_0$,  we make the following assumptions on $u_t$ to update the policy.


	\begin{assum}\label{ass:PE}
		For $t \geq t_0$, the input matrix $U_{0,t}$ is $\gamma {\sqrt{t(n+1)}}$-persistently exciting of order $1+n$, i.e., $\underline{\sigma}(\mathcal{H}_{n+1}(U_{0,t})) \geq \gamma {\sqrt{t(n+1)}}$, where $\gamma$ is a positive constant, and $\mathcal{H}_{n+1}(U_{0,t})$ is a Hankel matrix, i.e.,
		$$
		\mathcal{H}_{n+1}(U_{0,t}):=\begin{bmatrix} 
		u_0 & u_1 & \cdots & u_{t-n-1} \\
		u_1 & u_2 & \cdots & u_{t-n} \\
		\vdots & \vdots & \ddots & \vdots \\
		u_n & u_{n+1} & \cdots & u_{t-1}
		\end{bmatrix}.
		$$
	\end{assum}
	
	Assumption \ref{ass:PE} adopts a quantitative notion of PE~\cite[Definition 2]{coulson2022quantitative} and implies the rank condition \eqref{equ:rank}. {Here, the constant $\gamma$ is used to quantify the PE level. Since $\mathcal{H}_{n+1}(U_{0,t})\mathcal{H}_{n+1}^{\top}(U_{0,t})/t$ is the sample covariance of the input, the $\mathcal{O}(\sqrt{t})$ scaling in $\underline{\sigma}(\mathcal{H}_{n+1}(U_{0,t}))$ implies constant covariance excitation. Note that the PE assumption is universal in adaptive control~\cite{bian2016value, lewis2009reinforcement}.}
	\begin{assum}\label{ass:bound}
		There exist constants $\bar{u} > 0, \bar{x} > 0$ such that $ \|u_t\|\leq \bar{u}$ and $\|x_t\| \leq \bar{x} ,\forall t\in  \mathbb{N}$.
	\end{assum}

	{Assumption \ref{ass:bound} imposes the uniform boundedness and sequential stability \cite{cohen2018online} of the closed-loop system, i.e., the state does not blow up under the switching policy sequence. }
	 
	\begin{remark}\label{remark:assump23}
	{We use Assumptions \ref{ass:PE} and \ref{ass:bound} to decouple the convergence analysis and the input design problem. The latter concerns  designing  control inputs ensuring the excitation and boundedness assumptions, which is not the focus of this paper. Briefly, the input can be selected as
	\begin{equation}\label{equ:ut}
	u_t = K_t x_t + v_t, ~t \geq t_0,
	\end{equation}
	where $\{v_t\}$ is a probing noise sequence. In this case, Assumption \ref{ass:PE} holds with high probability when $\{v_t\}$ is i.i.d. Gaussian. Alternatively, a simple linear algebraic manipulation can also determine $v_t$, as the feedback term $K_tx_t$ is known. Of course, either approach will degrade the closed-loop performance. Assumption \ref{ass:bound} can possibly be circumnavigated by a more sophisticated sequential stability analysis \cite{cohen2018online}. Intuitively, proving sequential stability for \eqref{equ:ut} requires that (a) the LQR cost of the closed-loop system $A+BK_t$ is uniformly bounded, which can be achieved by our gradient methods; and (b) $K_t$ changes sufficiently slowly, which can be easily achieved by selecting a small stepsize $\eta$. We leave a detailed investigation of the consequences of these assumptions and the input design problem to future work, and instead focus on the convergence analysis in this paper.  \qed}
	\end{remark}


		

	Finally, we assume that the process noise $w_t$ is bounded as in \cite{de2019formulas,de2021low, dorfler2021certainty,dorfler22on, van2020noisy}, {which does not necessarily follow any particular statistics and can even be adversarial and correlated. }
	\begin{assum}\label{ass:noise}
		The process noise $w_t$ is upper-bounded, i.e., $\|w_t\| \leq \delta$ for some constant $\delta \geq 0$.
	\end{assum}
	
	{To better interpret our analysis and main results, we refer to $\gamma/\delta$ as the signal-to-noise ratio (SNR) describing the ratio between the useful and useless information $D_{0,t}$ and $W_{0,t}$. Such a notion of SNR is slightly different from the commonly used power-based definition.}
	
	While Algorithm \ref{alg:deepo} requires to compute the sample covariance matrices $\overline{X}_{0,t},\overline{U}_{0,t}, \overline{X}_{1,t}$ and $\Phi_{t}^{-1}$, it does not need to store all the historical data $(X_{0,t}, U_{0,t}, X_{1,t})$ and can be implemented recursively, as shown in the next subsection.

	\subsection{Recursive implementation of Algorithm \ref{alg:deepo}}\label{subsec:implement}
	We show how to efficiently implement Algorithm \ref{alg:deepo} via recursive rank-one update. First, the sample covariance matrices $\overline{X}_{0,t},\overline{U}_{0,t}, \overline{X}_{1,t}$ are updated recursively. To elaborate it, let
	$
	\psi_t := [u_{t}^{\top}, x_{t}^{\top}]^{\top}.
	$
	Then, $\overline{X}_{0,t+1}$ can be written as 
	$$
	\overline{X}_{0,t+1} = \frac{1}{t+1} 
	\sum_{k=0}^{t}x_k\psi_k^{\top} = \frac{t}{t+1}\overline{X}_{0,t} + \frac{1}{t+1}x_t\psi_t^{\top},
	$$
	where the first term is the weighted covariance matrix from the last iteration, and the second is a rank-one matrix. The other two matrices $\overline{U}_{0,t}, \overline{X}_{1,t}$ can be updated accordingly.
	
	Second, the covariance parameterization in (\ref{equ:k2v}) can also be implemented via rank-one update. We write the sample covariance  recursively  as
	$
	\Phi_{t+1} = (t\Phi_{t} + \psi_t\psi_t^{\top})/(t+1).
	$
	By Sherman-Morrison formula~\cite{sherman1950adjustment}, its inverse $\Phi_{t+1}^{-1}$ satisfies
	$$
	\Phi_{t+1}^{-1} = \frac{t+1}{t}\left(\Phi_{t}^{-1} - \frac{\Phi_{t}^{-1}\psi_t\psi_t^{\top}\Phi_{t}^{-1}}{t+\psi_t^{\top}\Phi_{t}^{-1}\psi_t}\right).
	$$
	Then, (\ref{equ:k2v}) can be written as a rank-one update
	\begin{align*}
	V_{t+1}&= \frac{t+1}{t}\left(\Phi_{t}^{-1} - \frac{\Phi_{t}^{-1}\psi_t\psi_t^{\top}\Phi_{t}^{-1}}{t+\psi_t^{\top}\Phi_{t}^{-1}\psi_t}\right) \Phi_{t} V_{t}' \\
	&= \frac{t+1}{t} \left(V_{t}' - \frac{\Phi_{t}^{-1}\psi_t\psi_t^{\top}V_{t}'}{t+\psi_t^{\top}\Phi_{t}^{-1}\psi_t}\right),
	\end{align*}
	where $\Phi_{t}^{-1}$ and $V_{t}'$ are given from the last iteration.

	\subsection{Non-asymptotic convergence guarantees for Algorithm \ref{alg:deepo}}\label{subsec:non}
	Let Assumptions \ref{ass:PE}-\ref{ass:noise} hold throughout this section. Let $T:=t-t_0+1$ denote the running time of Algorithm \ref{alg:deepo} and {$C^*:=\min_{K}C(K)$}. {We use the optimality gap $C(K_t)-C^*$ to quantify the convergence of $K_t$, which is aligned with the RL perspective \cite{fazel2018global,malik2019derivative,mohammadi2022convergence} and broader online optimization literature \cite{hazan2016introduction}. Define the average regret as the time-average optimality gap of the policy sequence}
	$$
	\text{Regret}_T: = \frac{1}{T}\sum_{t=t_0}^{t_0+T-1} (C(K_t)-C^*),~ T\geq 1.
	$$
	{The regret depends on $\{w_t\}$ but is not a well-defined random variable, since we do not assume any particular statistics on the noise. Now, we provide an upper bound on the regret for any bounded $\{w_t\}$ consistent with Assumption 4. }
	
	\begin{theorem}[\textbf{Non-asymptotic guarantees}]\label{thm:final}
	Given offline data $(X_{0,t_0}, U_{0,t_0}, X_{1,t_0})$, there exist constants $\nu_i>0,i\in\{1,2,3,4\}$  polynomial in $(\bar{x},\bar{u},\zeta,\omega,\|R\|,\underline{\sigma}(Q),\underline{\sigma}(R),C^*)$, such that, if $\eta \in (0, \nu_1]$ and $\text{SNR} \geq \nu_2$, then 
	$$
	\text{Regret}_T \leq \nu_3/{\sqrt{T}} + \nu_4  \text{SNR}^{-1/2}.
	$$ 
	\end{theorem}

	{Theorem \ref{thm:final} implies that, in the worst case, the time-average optimality gap converges sublinearly to a neighborhood of zero, the size of which scales inversely with the SNR. If the SNR tends to infinity with time, then DeePO is able to learn the optimal LQR gain.}

	
	{The proof of Theorem \ref{thm:final} is distinguished from the literature and much more challenging: (a) While we adopt the gradient descent framework as in zeroth-order PO \cite{fazel2018global,malik2019derivative,mohammadi2022convergence}, our cost function $J_t(V)$ is non-convex and time-varying due to the adaptive setting. Thus, we need to establish \textit{time-uniform projected} gradient dominance and smoothness properties. To prove convergence, we also need to show time-uniform boundedness of $J_t(V)$, which is not required in zeroth-order PO. (b) The convergence analysis for the gain sequence of indirect methods is mainly based on perturbation theory of Riccati equations~\cite{lu2023almost,mania_certainty_2019}. In comparison, we need to additionally characterize the gap between $J_t(V)$ and the actual cost $C(K_t)$ using perturbation analysis for Lyapunov equations. Moreover, the characterization of the optimality of $J_t^*$ requires Lemma 1, which is one of our key novel results. The main proof steps for Theorem \ref{thm:final} are sketched below.}

%
%

	\textbf{1) Quantifying the progress of the projected gradient descent in (\ref{equ:pgdV}) with a recursive form.} In the considered time-varying setting, the gradient dominance constant in Lemma \ref{lem:pl} involves a time-varying term $\underline{\sigma}(\overline{U}_{0,t})$. We provide its lower bound by using the robust fundamental lemma in \cite[Theorem 5]{coulson2022quantitative}, which quantifies $\underline{\sigma}(D_{0,t})$ under a sufficiently large SNR. 
	
	
	
	\begin{lemma}\label{lem:Ut}
		There exist positive constants $\zeta,\omega$ such that, 
		if
		$
		{\delta}/{\gamma } \leq \zeta/{(2\omega)}
		$,
		then $\underline{\sigma}(D_{0,t}) \geq \sqrt{t}\gamma \zeta/2$ and $\underline{\sigma}(\overline{U}_{0,t}) \geq \gamma \zeta/2$. 
	\end{lemma}

	{Here, the constant $\zeta$ is the system-theoretic parameter $\rho$ in \cite{coulson2022quantitative} that quantifies controllability, and $\omega$ is the parameter $\gamma$ in \cite{coulson2022quantitative} that quantifies effects of noise propagation.} By Lemmas \ref{lem:pl} and \ref{lem:Ut}, it is straightforward to show that $J_t(V)$ has a uniform gradient dominance constant  independent of $t$.
	
	\begin{lemma}[\textbf{Uniform projected gradient dominance}]\label{lem:uni_pl}
		If $V \in \mathcal{S}_{t}(a)$ and
		$
		{\delta}/{\gamma } \leq  {\sqrt{n+1}}\zeta/{(2n\omega)}$,
		then it holds that 
		$$
		J_t(V) - J_t^* \leq \mu(a) \|\Pi_{\overline{X}_{0,t}} \nabla J_t(V)\|,
		$$
		where   $\mu(a)={8a^2}/({\gamma \zeta(\underline{\sigma}(Q))^{{3}/{2}}(\underline{\sigma}(R))^{{1}/{2}}}).$
	\end{lemma}
	
	Next, we show that $J_t(V)$ also has a uniform smoothness constant, which follows directly from Lemma \ref{lem:smooth} and the boundedness of the data matrices $\|\overline{U}_{0,t}\| \leq  \bar{u}(\bar{u}+\bar{x})$ and  $\|\overline{X}_{1,t}\| \leq  \bar{x}(\bar{u}+\bar{x})$ in Assumption \ref{ass:bound}.
	
	\begin{lemma}[\textbf{Uniform local smoothness}]\label{lem:smooth_V}
		For any $V,V' \in \mathcal{S}_t(a)$ satisfying $V+\phi(V'-V) \in \mathcal{S}_t(a)$ for all scalar $\phi \in [0,1]$, it holds that
		$
		J_t(V') \leq J_t(V) +  \langle \nabla J_t(V), V' - V \rangle + {l(a)}\|V'-V\|^2/2,
		$
		where the uniform smoothness constant is
		\begin{align*}
		&l(a) = 2\left(\bar{u}^2 + \bar{x}\bar{u}\right)^2\frac{a\|R\|}{\underline{\sigma}(Q)} +\left(\bar{x}^2 + \bar{x}\bar{u}\right)^2 \frac{2a^2}{\underline{\sigma}^2(Q)} \\
		&~~~\times \left(2(\bar{u}^2 + \bar{x}\bar{u})^2\|R\| + 2(\bar{x}^2 + \bar{x}\bar{u})^2a + 2a -\underline{\sigma}(Q)\right).
		\end{align*}
	\end{lemma}
	
	Following the proof of Theorem \ref{thm:conv},  for $\eta \in [0, 1/l(J_t(V_t))]$, the progress of (\ref{equ:pgdV}) can be characterized as 
	\begin{equation}\label{equ:M3}
	\begin{aligned}
	&J_t(V_t')-J_t(V_t) \\
	& \leq  - \frac{\eta}{\mu^2(J_t(V_t))}\left(1 -\frac{l(J_t(V_t))\eta}{2}\right)(J_t(V_t) - J_t^*)^2.
	\end{aligned}
	\end{equation}
	
	
	\textbf{2) Bounding the difference among $C(K_t), J_{t+1}(V_{t+1})$ and $J_t(V_t')$}. 
	Since the noise matrix is disregarded in the parameterized closed-loop matrix \eqref{equ:clos}, the parameterized LQR cost $J_{t+1}(V_{t+1})$ and $J_t(V_t')$ are not equal to $C(K_t)$. To bound their difference, we quantify the distance between their closed-loop matrices. In particular, it holds that
	\begin{equation}\label{equ:m1m2}
	\begin{aligned}
	&\left\| A+BK_t -\overline{X}_{1,t}V_t' \right\|\\
	&\overset{\eqref{equ:dynamics}}{=} \left\|   A+BK_t -\left(\begin{bmatrix}
	B & A
	\end{bmatrix} \Phi_t + \overline{W}_{0,t}\right)V_t' \right\|\\
	&\overset{\eqref{equ:KV'}}{=} \left\| {W}_{0,t}D_{0,t}^{\dagger}\begin{bmatrix}
	K_{t} \\ I_n
	\end{bmatrix} \right\|.
	\end{aligned}
	\end{equation}
	By Lemma \ref{lem:Ut} and Assumption \ref{ass:noise}, it follows that
	\begin{equation}\label{equ:snr}
	\left\|W_{0,t}D_{0,t}^{\dagger}\right\| \leq  {\|W_{0,t}\|}{(\underline{\sigma}(D_{0,t}))^{-1}} \leq \frac{2\delta}{\gamma \zeta}, 
	\end{equation}
	and (\ref{equ:m1m2}) further leads to
	\begin{equation}\label{equ:diffff}
	\left\| A+BK_t -\overline{X}_{1,t}V_t' \right\| \leq 2(1+\|K_{t}\|){\delta}/{(\gamma \zeta)}.
	\end{equation}
	Then, we can use perturbation theory of Lyapunov equations to bound $|C(K_t)-J_t(V_t')|$. Let $p(a)$ be the   polynomial of a scalar $a\geq 0$:
	$$
	p(a) = \frac{16a^3}{\underline{\sigma}^2(Q)\zeta}\left(1+\frac{a}{\underline{\sigma}(Q)}\right) \left(1+\sqrt{\frac{a}{\underline{\sigma}(R)}}\right).
	$$
	Then, we have the following results.
	
	\begin{lemma}\label{lem:cost_diff}
		For $t\geq t_0$, if $V_t' \in \mathcal{S}_t$ and
		$$
		\frac{\delta}{\gamma} \leq \min\left\{\frac{J_t^2(V_t') \zeta}{\underline{\sigma}(Q)p(J_t(V_t'))},\frac{ \zeta}{2\omega}\right\},
		$$
		then it holds that  
		\begin{align}
		&|C(K_{t}) - J_{t}(V_t')| \leq  \frac{ p(J_t(V_t'))\delta}{2\gamma }, \label{equ:1}\\
		&|J_{t+1}(V_{t+1})- J_{t}(V_t')| \leq  \frac{ p(J_t(V_t'))\delta}{	\gamma },\label{equ:2}\\
		&|J_{t+1}(V_{t+1}) - C(K_{t})| \leq  \frac{3p(J_t(V_t'))\delta}{2	\gamma }.\label{equ:3}
		\end{align}
	\end{lemma}
	
	\textbf{3) Bounding $|J_t^*-C^*|$}. By Lemma \ref{lem:equi}, it is equivalent to provide an upper bound for $|C_{\text{CE},t}^* - C^*|$, which is exactly the sub-optimality gap of the certainty-equivalence LQR \eqref{prob:indirect}. Since by (\ref{equ:snr})  $\|[\widehat{B}_t,\widehat{A}_t]-[B,A]\| = \|W_{0,t}D_{0,t}^{\dagger}\| \leq 2\delta/(\gamma \zeta)$, we can use perturbation theory of Riccati equations~\cite[Proposition 1]{mania_certainty_2019} to provide an upper bound for $|C_{\text{CE},t}^* - C^*|$.
	
	\begin{lemma}\label{lem:opt_diff}
		There exist two constants $c_1,c_2$ depending only on $n, A,B,Q,R$ such that, if $\delta/\gamma  \leq \min\{c_1\zeta, \zeta/({2\omega})\}$, then
		$
		\left|J_t^* -C^* \right| \leq c_2\delta/(\gamma \zeta).
		$
	\end{lemma}
	
	\textbf{4) Uniform bounds for $J_t(V)$}. 
	Both $\mu(J_t(V_t)),l(J_t(V_t))$ in (\ref{equ:M3}) and the bounds in Lemma \ref{lem:cost_diff} depend on $J_t(V_t)$ and $J_t(V_t')$. Next, we provide uniform bounds for $J_t(V_t)$ and $J_t(V_t')$. For brevity, let $c_3 :=  l(\underline{\sigma}(Q)), ~c_4:= \mu(\underline{\sigma}(Q)),
	\bar{J} :=J_{t_0}^*+ C^* + c_1c_2 + {1}/{(2c_3c_4^2)}+ 1,
	\bar{\mu}: = \mu(\bar{J}), \bar{l}: = l(\bar{J}), \text{and}~ \bar{p}: = p(\bar{J}).$
	Then, we have the following result.

	\begin{lemma}\label{lem:boundcost}
		If $\eta \in (0, 1/\bar{l}]$ and
		\begin{equation}\label{equ:cond_bound}
		\frac{\delta}{\gamma } \leq \min\left\{
		\frac{\bar{J}^2}{\underline{\sigma}(Q)\bar{p}},\frac{\zeta}{2\omega},\frac{\eta}{2\bar{\mu}^2\bar{p}}, c_1\zeta \right\},
		\end{equation}
		then $J_t(V_t)\leq \bar{J}$ for $t > t_0$ and $J_t(V_t')\leq \bar{J}$ for $t \geq t_0$.
	\end{lemma}
	
	Equipped with Lemmas \ref{lem:Ut}-\ref{lem:boundcost}, Theorem \ref{thm:final} can be proved combining standard gradient descent arguments~\cite{hazan2016introduction}.
	
	\subsection{Discussion}\label{subsec:discuss} 
	The DeePO method in Algorithm \ref{alg:deepo} is adaptive and direct  in that it learns from online closed-loop data without any explicit SysID. This is in contrast to the indirect adaptive control that involves an identification step~\cite{ campi1998adaptive, wang2021exact, lu2023almost}, and the episodic methods using single or multiple alternating episodes of data collection and control~\cite{abbasi2011regret,fazel2018global,de2019formulas}. 
	
	
	The DeePO method has a recursive policy update and is computationally efficient. In particular, Algorithm \ref{alg:deepo} performs only one-step projected gradient descent per time efficiently using online closed-loop data, and can be implemented recursively (see Section \ref{subsec:implement}). {In comparison, the indirect adaptive control~\cite{ campi1998adaptive, wang2021exact, lu2023almost} requires to solve one certainty-equivalence LQR problem (\ref{prob:indirect}) per time. While one can do one or a few Riccati iterations per
	time step to gain in terms of computation time, it remains unclear if such an iterative approach has provable theoretical guarantees.}

	
	Our analysis of Theorem \ref{thm:final} is non-asymptotic and independent of the noise statistics. It improves over single batch methods~\cite{de2021low,dorfler2021certainty,dorfler22on}, where their performance also depends on SNR but does not decay over time \cite[Theorem 4.1]{dorfler22on}. Our sublinear decrease matches the attainable rate $\mathcal{O}(1/\sqrt
	{T})$ of first-order methods in online convex optimization of smooth functions~\cite[Chapter 3]{hazan2016introduction}, even though our considered LQR problem (\ref{equ:transfer})-(\ref{equ:Sigma}) is non-convex~\cite{fazel2018global}. This indicates a favorable sample efficiency of DeePO to learn from online closed-loop data. {The regret bound also has a polynomial dependence on the optimal LQR cost $C^*$, which is a key system-theoretic parameter\cite{tsiamis2023statistical}. This is in line with the scaling in terms of $C^*$ of the regret bounds in indirect adaptive control~\cite{simchowitz2020naive}.}

	
	
	{While Theorem \ref{thm:final} assumes the boundedness of noise, DeePO can handle stochastic noise as well. Suppose that $\{w_t\}$ is white Gaussian, then the SNR has an explicit high-confidence bound, i.e., it decays inversely proportional to the square root of time [4]. Together with Theorem 2, it can be shown that the regret converges sublinearly with high probability. A detailed investigation of the stochastic case is left to future work.}

	{Since the policy update of DeePO is end-to-end, we can enhance the adaptivity by using sliding window data or forgetting factor, which is useful in time-varying systems. For example, adding a forgetting factor $\beta \in (0,1)$ to the sample covariance \eqref{def:cov} leads to $\Phi:= D_0 S D_0^{\top}/t$ with $S=\text{diag}(\beta^{t-1},\beta^{t-2},\cdots,1)$, and all results in Section \ref{subsec:non} continue to hold. We defer a rigorous analysis of time-varying systems to future work.}
	

		

	\section{Simulations}\label{sec:simu}
	In this section, we first simulate over a randomly generated linear system to validate the convergence of DeePO for learning the LQR using offline and online closed-loop data, respectively. Then, we compare DeePO with indirect adaptive control \cite{wang2021exact, lu2023almost} and zeroth-order PO \cite{fazel2018global,malik2019derivative,mohammadi2022convergence}. 
	\subsection{Convergence of DeePO for the LQR with offline data}
\begin{figure}[t]
		\centerline{\includegraphics[width=55mm]{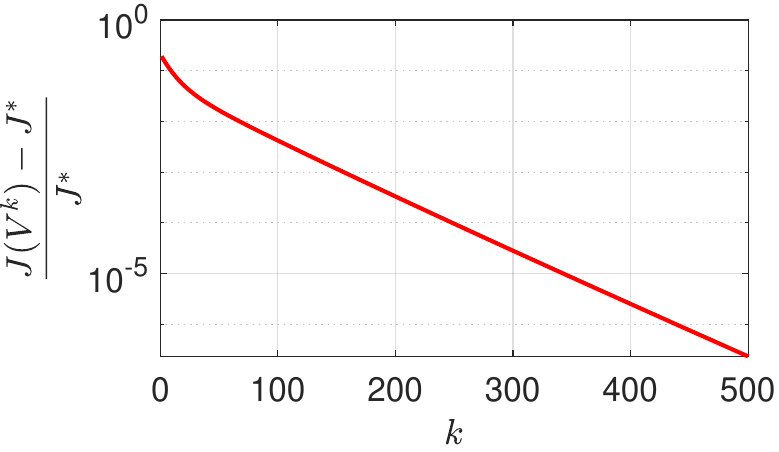}}
		\caption{Convergence of DeePO for the LQR with offline data, where $J^*$ is obtained by solving the convex program of  (\ref{prob:equiV}).}
		\label{pic:conv}
	\end{figure}

	We randomly generate a controllable and open-loop stable system $(A,B)$ with $n=4,m=2$ as
	\begin{align*}
	&A = \begin{bmatrix}
	-0.13 &  0.14  &  -0.29  & 0.28\\
	0.48  & 0.09 &  0.41 &  0.30\\
	-0.01  &  0.04  &  0.17  & 0.43\\
	0.14  &  0.31 &  -0.29 &  -0.10
	\end{bmatrix}, B = \begin{bmatrix}
	1.63  &  0.93\\
	0.26 &  1.79\\
	1.46  & 1.18\\
	0.77  & 0.11
	\end{bmatrix}.
	\end{align*}
	{Let $Q = I_4$ and $R = I_2$.} We generate PE data $X_{0},{U}_{0},W_{0}$ of length $8$ from a standard normal distribution and compute $X_{1}$ using the linear dynamics (\ref{equ:dynamics}). {The SNR (computed as $\underline{\sigma}(D_{0})/\|W_{0}\|$) is $-0.12$ dB.} We use only $(X_{0},{U}_{0},X_{1})$ to perform DeePO (\ref{equ:gd}) to solve the LQR problem with covariance parameterization (\ref{prob:equiV}), which is feasible under the given set of PE data. We set the stepsize  to $\eta = 0.1$ and the initial policy to
	$
	V^0 = \Phi^{-1}[0,I_4]^{\top}
	\in \mathcal{S}.
	$ Fig. \ref{pic:conv} shows that the LQR cost converges to $J^*$ at a linear rate, which implies that the sublinear rate certified in Theorem \ref{thm:conv} may be conservative. 
	
	\subsection{Convergence of DeePO for the adaptive learning of the LQR with online closed-loop data}\label{subsec:simu_online}
	\begin{figure}[t]
		\centerline{\includegraphics[width=55mm]{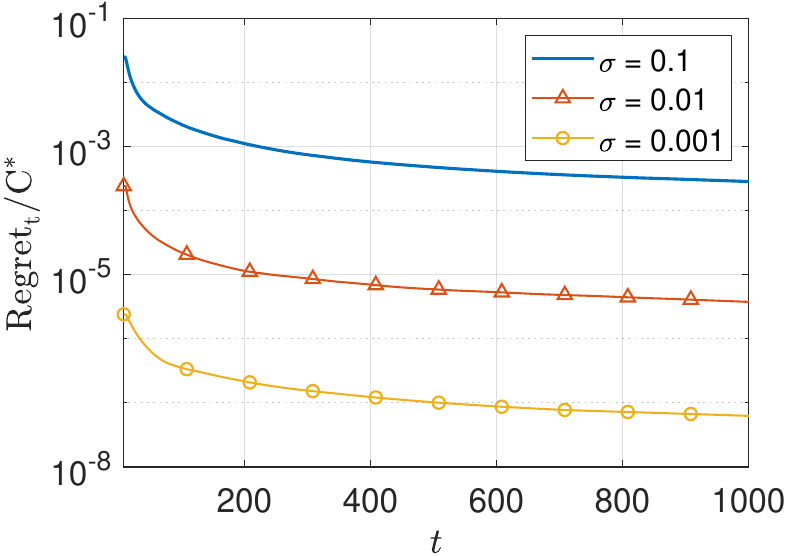}}
		\caption{Convergence of DeePO for adaptive learning of the LQR under different level of noise $\sigma$.}
		\label{pic:adap}
	\end{figure}
	In this subsection, we perform Algorithm \ref{alg:deepo} for the adaptive learning of the LQR with online closed-loop data to validate   Theorem \ref{thm:final}. We set $x_0 = 0$ and $u_t\sim \mathcal{N}(0,I_2)$ for $t<t_0$ with $t_0=8$. For $t \geq t_0$, we set $u_t =K_tx_t + v_t$ with $v_t \sim \mathcal{N}(0,I_2)$ to ensure persistency of excitation. The white noise sequence $\{w_t\}$ is drawn from a uniform distribution, where each element of $w_t$ is uniformly sampled from $[0, \sigma]$. {We consider three  noise levels $\sigma \in \{0.1,0.01,0.001\}$, which correspond approximately to the SNR $\in [0,5],[10,15], [20,25]$ dB (computed as $\underline{\sigma}(D_{0,t})/\|W_{0,t}\|$), respectively. Note that the SNR varies within a interval due to the random online sample of noise. } We set the stepsize to $\eta = 0.01$.

	\begin{figure}[t]
		\centerline{\includegraphics[width=55mm]{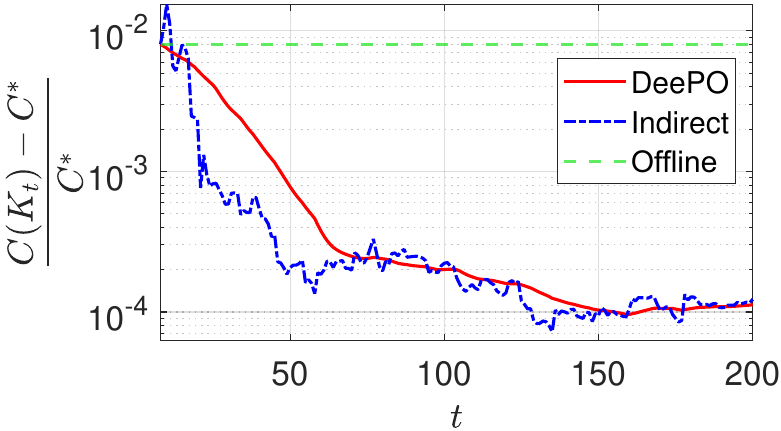}}
		\caption{Comparison of DeePO and indirect methods for adaptive learning of the LQR. We also plot the optimality gap of the initial gain (dashed green line) computed by single batch methods from offline data~\cite{de2021low,dorfler2021certainty,dorfler22on}.}
		\label{pic:comp}
	\end{figure}
	
	Fig. \ref{pic:adap} shows that average regret matches the expected sublinear decrease $\mathcal{O}(1/\sqrt{T})$ in Theorem \ref{thm:final}. Moreover, since the bounded noise $w_t$ has non-zero mean, the regret does not converge to zero. Roughly, the bias scales as SNR$^{-2}$ rather than the expected SNR$^{-1/2}$, which implies that our bound on the bias term certified in Theorem \ref{thm:final} may be conservative. Bridging this gap requires tighter perturbation bounds.

	\subsection{Comparison with indirect certainty-equivalence adaptive control for adaptive learning of the LQR}\label{subsec:comp_ind}
	Consider the system proposed in \cite[Section 6]{dean2020sample}
	\begin{equation}\label{equ:benchmark}
	\begin{aligned}
	&A = \begin{bmatrix}
	1.01 &  0.01  &  0\\
	0.01  & 1.01 &  0.01\\
	0  &  0.01  &  1.01
	\end{bmatrix}, ~B = I_3,
	\end{aligned}
	\end{equation}
	which corresponds to a discrete-time marginally unstable Laplacian system. Let $Q = R = I_3$. {To verify that DeePO can handle unbounded noise, let $w_t\sim  \mathcal{N}(0, I_3/100)$ which corresponds to SNR $\in [0,5]$ as in \cite[Section VI]{dorfler2021certainty}}. We compare Algorithm \ref{alg:deepo} with the indirect adaptive approach in~\cite{wang2021exact, lu2023almost}. {Specifically, the indirect adaptive method alternates between finding the certainty-equivalence LQR gain $K_t$ \eqref{prob:indirect} every time step and using $u_t =K_tx_t + v_t$  with $v_t \sim \mathcal{N}(0,I_2)$ to obtain the new state $x_{t+1}$.} For both methods, we use the solution to the LQR problem with covariance parameterization \eqref{prob:equiV} based on $(X_{0,t_0}, U_{0,t_0}, X_{1,t_0})$ as the initial gain $K_{t_0}$.

	\begin{figure}[t]
		\centerline{\includegraphics[width=55mm]{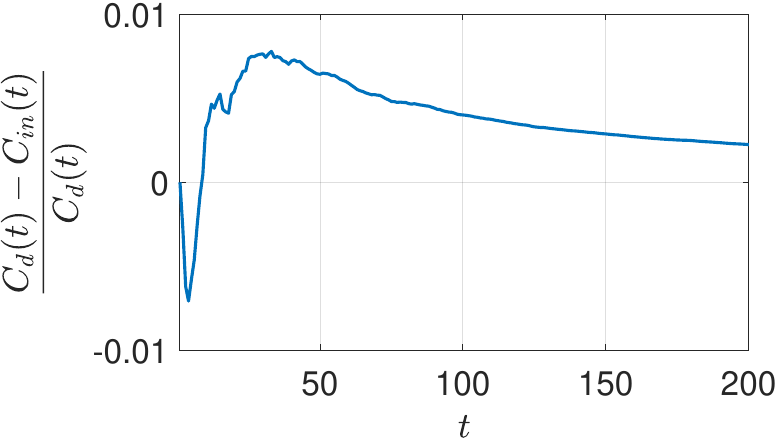}}
		\caption{{Comparison of the finite-horizon cost of DeePO and the indirect method.}}
		\label{pic:finite}
	\end{figure}

	Fig. \ref{pic:comp} shows that the optimality gap of both methods converges sublinearly to $10^{-4}$ in $200$ time steps. By using online closed-loop data, both methods improve the performance over the initial gain. For $t \leq20$, the indirect approach exhibits faster convergence than DeePO. This is because DeePO  only performs one-step projected gradient descent per time, and hence requires more iterations to asymptotically approach the certainty-equivalence LQR. Indeed, they achieve similar performance for $t\geq 60$. Moreover, thanks to the recursive policy update, the curve of DeePO is significantly smoother.
	
	\begin{figure}[t]
		\centerline{\includegraphics[width=55mm]{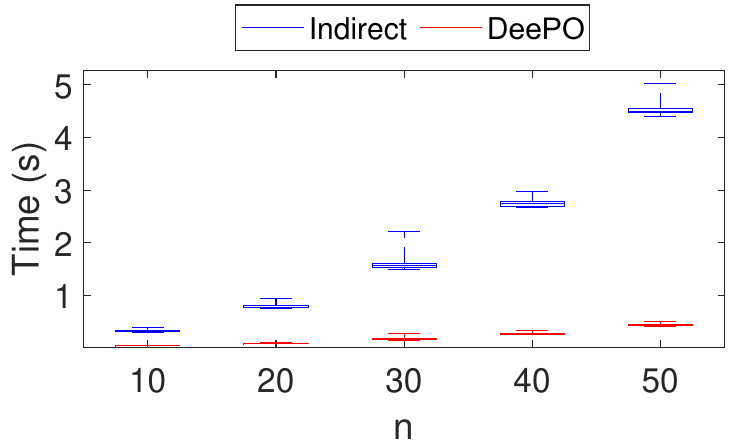}}
		\caption{Comparison of the computation time for performing $100$ time steps.}
		\label{pic:time_order}
	\end{figure}


	\begin{figure}[t]
		\centerline{\includegraphics[width=55mm]{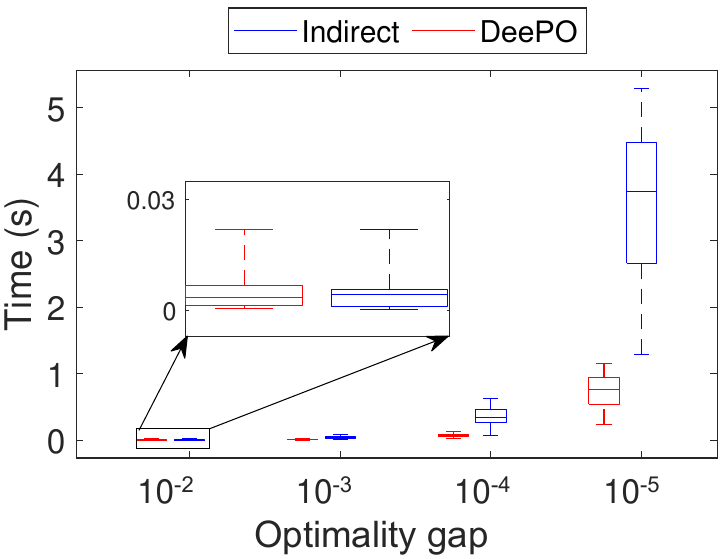}}
		\caption{Comparison of the computation time to achieve a specified optimality gap.}
		\label{pic:acc}
	\end{figure}
	
	{Next, we compare the finite-horizon cost $\sum_{k=0}^{t-1}\|z_k\|^2$ induced by DeePO and the indirect method, denoted by $C_d(t)$ and $C_{in}(t)$, respectively. Fig. \ref{pic:finite} reports the average finite-horizon cost from $50$ independent trials, showing that their performance is extremely close during the entire run time. Further, at the initial several steps, the DeePO method has slightly lower cost due to smoother switching of the policy. After that, the indirect method performs better but their relative cost approaches to zero asymptotically.}

	Finally, we compare their efficiency in terms of real-time computation.
	For a fair comparison, we apply rank-one update for Algorithm \ref{alg:deepo} (see Section \ref{subsec:implement}) and recursive least
	squares for  SysID (\ref{equ:sysid}). {Let $Q$ and $R$ be identity matrices with proper dimension.} With identified model, the certainty-equivalence LQR problem \eqref{prob:indirect} is solved using the \textit{dlqr} function in MATLAB. First, we evaluate the computation time for systems with different dimensions during $100$ time steps. {Let the dimension of the state and input be equal, and $n =m \in \{10,20,30,40,50\}$. } For each instance, we perform $50$ independent trials, where each trial uses randomly generated stable state matrix $A$ and identity input matrix $B$. {The box plot in Fig. \ref{pic:time_order} shows that DeePO is significantly more efficient than indirect methods, and their gap scales significantly with the dimension.} Second, we also compare the computation time to achieve different optimality gap $\epsilon := {(C(K)-C^*)}/C^*$ for a fixed system dimension $n = 4$. Fig. \ref{pic:acc} shows the results from $50$ independent trials. For a large optimality gap (e.g., $\epsilon = 10^{-2}$), the computation time is similar for both methods. This reveals that while the indirect method converges faster at the beginning (see Fig. \ref{pic:comp}), its real-time computation is more time-consuming. For a small optimality gap (e.g., $\epsilon=10^{-5}$), DeePO requires significantly less computation time. Since a rank-one update has been applied in both methods, the difference of computational efficiency is due to DeePO only performing one-step gradient descent per time, but the indirect method needs to solve a Riccati equation.

	\subsection{Comparison with zeroth-order PO}
	\begin{table}[t]\label{table}
		\begin{center}
			\begin{tabular}{|c|c|c|c|}
				\hline
				& $\epsilon = 1$ & $\epsilon =0.1$ & $\epsilon =0.01$ \\
				\hline
				Zeroth-order PO (\# of trajectories) & $1393$ & $45260$ & $151607$\\
				\hline
				DeePO (\# of input-state pairs) & $10$ & $25$ & $49$\\
				\hline
			\end{tabular}
		\end{center}
		\caption{Comparison of the sample complexity to achieve a specified optimality gap $\epsilon$.}
	\end{table}
	
	While the  zeroth-order PO methods\cite{fazel2018global,malik2019derivative,mohammadi2022convergence} are episodic (see Fig. \ref{subfig:episodic}), we compare them with DeePO in Algorithm \ref{alg:deepo} to demonstrate our sample efficiency. {We use the simulation model (\ref{equ:benchmark}) with $Q = 10\times  I_3$, $R = I_3$ and  $w_t \sim \mathcal{N}(0, 0.01I_n)$.} We adopt the zeroth-order PO method in (\ref{equ:pogd})-(\ref{equ:gdest}) for comparison. In particular, we use a minibatch of $30$ zeroth-order samples for (\ref{equ:gdest}) to reduce the variance of gradient estimate. We set the smooth radius to $r = 0.02$, the stepsize to $\eta = 10^{-3}$, and the length of the trajectory to $T=50$. The setting for Algorithm \ref{alg:deepo} is the same as in Section \ref{subsec:simu_online}. For both methods, we use $-0.15\times I_3$ as the initial stabilizing gain. 
	
	Table I demonstrates the sample complexity of zeroth-order PO (in terms of number of trajectories) and DeePO (in terms of number of input-state pairs) to achieve different optimality gap. It indicates that the zeroth-order PO is vastly less efficient for solving the LQR problem, which is in line with the sample complexity discussion in \cite{tu2019gap}.
	
	\section{Concluding remarks}\label{sec:conclution}
	This paper proposed DeePO for the direct adaptive learning of the LQR based on the covariance parameterization. The proposed method is direct, adaptive, with closed-loop data, and has a recursive implementation. Hence, we provided a viable angle of attack to the open problem in \cite{markovsky2021behavioral,markovsky2023data,annaswamy2023adaptive}. 
	
	We believe that our paper leads to fruitful future works. {As discussed in Remark \ref{remark:assump23}, input design, quantifying the effects of probing noise, and weakening Assumptions \ref{ass:PE} and \ref{ass:bound} are formidable problems worthy of further investigations. We believe that the sequential stability analysis \cite{cohen2018online} would be the key to achieving this.} As observed in the simulation, the sublinear rate in Theorem \ref{thm:conv} and the dependence on the SNR in Theorem \ref{thm:final} may be conservative, and their sharper analysis is an important future work. {While this paper considers only the LQR problem, we believe that the covariance parameterization and DeePO can be extended for other objectives (e.g., output feedback control), performance indices (e.g., $\mathcal{H}_{\infty}$ norm), and system classes (e.g., time-varying systems). It is valuable to move beyond the bounded noise assumption and analyze the expected regret under stochastic noise.} The adaptivity of DeePO can be further enhanced by using sliding data window and/or forgetting factor, and a rigorous analysis is also open.
	

	%

	%

	\section*{Acknowledgement}
	We would like to thank Prof. Mihailo R. Jovanovi{\'c} from University of Southern California and Dr. Hesameddin Mohammadi from NVIDIA for fruitful discussions.

	\appendices
	\section{Proof of Lemma \ref{lem:equi}}
	By the covariance parameterization (\ref{equ:newpara}) and the positive definiteness of $\Phi$ (due to the PE condition \eqref{equ:rank}), $V$ is uniquely determined as $V = \Phi^{-1}
	[K^{\top},I_n]^{\top}$. Then, the LQR problem with covariance parameterization \eqref{prob:equiV} can be reformulated as 
	\begin{align*}
	&\mathop{\text{minimize}}\limits_{K, \Sigma\geq 0} ~\text{Tr}\left((Q+K^{\top}RK)\Sigma\right),\\
	&\text{subject to} ~~	\Sigma = I_n + X_1D_0^{\dagger}\begin{bmatrix}
	K \\
	I_n
	\end{bmatrix}\Sigma \begin{bmatrix}
	K \\
	I_n
	\end{bmatrix}^{\top}(X_1D_0^{\dagger})^{\top}.
	\end{align*}
	By the definition $[\widehat{B},\widehat{A}] = X_1D_0^{\dagger}$ in \eqref{equ:sysid}, this problem is exactly the certainty-equivalence LQR problem \eqref{prob:indirect}. 
	
	\section{Proof in Section \ref{sec:deepo}}
	We provide some useful bounds which will be invoked frequently in the rest of the paper.
	\begin{lemma}\label{lem:bounds}
		For $V \in \mathcal{S}$, it holds
		(\romannumeral1) 
		$
		\|\Sigma_V\|  \leq {J(V)}/{\underline{\sigma}(Q)},
		$
		(\romannumeral2) 
		$
		\|P_V\| \leq J(V),
		$
		and (\romannumeral3) 
		$
		\|\overline{U}_0V\|_F \leq ({J(V)}/{\underline{\sigma}(R)})^{{1}/{2}}.
		$
	\end{lemma}
	\begin{proof}
		The proof of (\romannumeral1) and (\romannumeral2) follows directly from the definition of $J(V)$. To show (\romannumeral3), we have
		$ J(V) = \text{Tr}((Q+V^{\top}\overline{U}_0^{\top}R\overline{U}_0V)\Sigma_V)
		\geq \text{Tr}(V^{\top}\overline{U}_0^{\top}R\overline{U}_0V\Sigma_V)
		\geq \text{Tr}(V^{\top}\overline{U}_0^{\top}\overline{U}_0V) \underline{\sigma}(R) 
		= \|\overline{U}_0V\|_F^2 \underline{\sigma}(R).$
	\end{proof}	
	
	\subsection{Proof of Lemma \ref{lem:gradient}}
	The proof follows the same vein as that of \cite[Proposition 1]{maartensson2009gradient}. We first compute the differential of $P_V$. Define $E_V := (\overline{U}_0^{\top}R\overline{U}_0+\overline{X}_1^{\top}P_V\overline{X}_1)V$. By definition, it holds that
	$
	dP_V =  dV^{\top} E_V  + E^{\top}_VdV + V^{\top}\overline{X}_1^{\top}dP_V\overline{X}_1V.
	$ 
	
	Since $\rho(\overline{X}_1V)<1$, we obtain
	\begin{align*}
	dP_V &=\sum_{i=0}^{\infty} \left(V^{\top}\overline{X}_1^{\top}\right)^i\left(dV^{\top} E_V  + E^{\top}_VdV \right)\left(\overline{X}_1V\right)^i\\
	\text{Tr}(dP_V) &= \text{Tr}\left(2\sum_{i=0}^{\infty} (\overline{X}_1V)^i (V^{\top}\overline{X}_1^{\top})^i\cdot E^{\top}_VdV\right) \\
	&=\text{Tr}\left(2\Sigma_V E_V^{\top} dV \right).
	\end{align*}
	
	\subsection{Proof of Lemma \ref{lem:pl}}
	We apply \cite[Theorem 1]{sun2021analysis} to prove Lemma \ref{lem:pl}. We first relate the LQR problem with covariance parameterization (\ref{prob:equiV}) to the following convex reparameterization via a change of variables $V=L\Sigma^{-1}$:
	\begin{equation}\label{prob:convex}
	\begin{aligned}
	&\mathop{\text{minimize}}\limits_{L,\Sigma} ~f(L,\Sigma):= \text{Tr}(Q\Sigma) + \text{Tr}(L\Sigma^{-1}L^{\top}\overline{U}_0^{\top}R\overline{U}_0), \\
	&\text{subject to} ~~\Sigma = \overline{X}_0L,~
	\begin{bmatrix}
	\Sigma-I_n & \overline{X}_1L \\
	L^{\top}\overline{X}_1^{\top} & \Sigma
	\end{bmatrix} \succeq 0.
	\end{aligned}
	\end{equation}	
	Let $\mathcal{L}$ denote the feasible set of (\ref{prob:convex}). The equivalence between (\ref{prob:equiV}) and (\ref{prob:convex}) are established below.
	
	\begin{lemma}\label{lem:asmp1}
		For any $(L,\Sigma) \in \mathcal{L}$, it holds that $\Sigma$ is invertible and $L\Sigma^{-1} \in \mathcal{S}$. Moreover, for $V\in \mathcal{S}$, we have
		\begin{equation}\label{equ:equi}
		J(V) = \min_{L,\Sigma} \{f(L,\Sigma),\text{s.t.} (L,\Sigma) \in \mathcal{L}, L\Sigma^{-1} = V\}.
		\end{equation}
	\end{lemma}
	\begin{proof}
		Applying the Schur complement to the linear matrix inequality constraint in \eqref{prob:convex} yields $\Sigma \succ 0$ and
		$
		\Sigma-I_n -\overline{X}_1L\Sigma^{-1}L^{\top}\overline{X}_1^{\top} \succeq 0.
		$
		By the non-singularity of $\Sigma$, let $V = L\Sigma^{-1}$. Then, a substitution of $L = V\Sigma$ into the above inequality yields
		$
		\Sigma-I_n -\overline{X}_1V\Sigma V^{\top}\overline{X}_1^{\top} \succeq 0.
		$
		Thus, $\overline{X}_1V$ is stable, i.e., $\rho(\overline{X}_1V)<1$. Since the first constraint of \eqref{prob:convex} implies $\overline{X}_0V = \overline{X}_0L\Sigma^{-1} = \Sigma \Sigma^{-1} = I_n$, it holds that $V = L\Sigma^{-1} \in \mathcal{S}$.
		
		Next, we prove the second statement. Using the constraint $V = L\Sigma^{-1}$ and the Schur complement, the right-hand side of (\ref{equ:equi}) becomes
		\begin{equation}\label{prob:JG}
		\begin{aligned}
		&\min_{\Sigma \succ 0} ~ \text{Tr}\left(\left(Q+V^{\top}\overline{U}_0^{\top}R\overline{U}_0V\right)\Sigma\right)\\
		&~\text{s.t.} ~\overline{X}_0V = I_n, \Sigma \succeq I_n + \overline{X}_1V\Sigma V^{\top}\overline{X}_1^{\top}.
		\end{aligned}
		\end{equation}
		Let $\Sigma(\Theta)$ be the unique positive definite solution of the Lyapunov equation
		$
		\Sigma(\Theta) = \Theta + \overline{X}_1V\Sigma(\Theta) V^{\top}\overline{X}_1^{\top}
		$
		with $\Theta \succeq I_n$. By monotonicity of $\Sigma(\Theta)$, we have $\Sigma(\Theta) \succeq \Sigma(I_n)$. Since $Q+V^{\top}\overline{U}_0^{\top}R\overline{U}_0V \succ 0$, the minimum of (\ref{prob:JG}) is attained at $\Sigma(I_n)$, which is
		$
		\text{Tr}((Q+V^{\top}\overline{U}_0^{\top}R\overline{U}_0V)\Sigma(I_n))
		$
		with $\overline{X}_0V = I_n$. This is exactly the definition of $J(V)$.
	\end{proof}

	Next, we show the convexity of the reparameterization \eqref{prob:convex}.
	{\begin{lemma}\label{lem:asmp2}
			The feasible set $\mathcal{L}$ is convex, and $f(L,\Sigma)$ is convex over $\mathcal{L}$. Moreover, $f(L,\Sigma)$ is differentiable over an open domain that contains $\mathcal{L}$. 
	\end{lemma}}
	\begin{proof}
		Since the constraints in (\ref{prob:convex}) are linear in $(L,\Sigma)$, the feasible set $\mathcal{L}$ is convex. Clearly,  $f(L,\Sigma)$ is differentiable over $\mathcal{L}$. Hence, $f$ is convex over $\mathcal{L}$. Define the derivative and Hessian of $f$ along the direction $(\tilde{L},\tilde{\Sigma})$ as
		\begin{align*}
		&f'(L,\Sigma)[(\tilde{L},\tilde{\Sigma})]:=\frac{d}{d\tau}f(L+\tau\tilde{L},\Sigma+\tau\tilde{\Sigma})\big|_{\tau=0},\\
		&\nabla^2 f(L,\Sigma)[(\tilde{L},\tilde{\Sigma}),(\tilde{L},\tilde{\Sigma})]:= \frac{d^2}{d\tau^2}f(L+\tau\tilde{L},\Sigma+\tau\tilde{\Sigma})\big|_{\tau=0},
		\end{align*}
		respectively. By the definition of $f(L,\Sigma)$, it holds that
		$
		f'(L,\Sigma)[(\tilde{L},\tilde{\Sigma})] = \text{Tr}(Q\tilde{\Sigma})
		- \text{Tr}(\Sigma^{-1}L^{\top}\overline{U}_0^{\top}R\overline{U}_0L\Sigma^{-1}\tilde{\Sigma}) 
		+ 2\text{Tr}(\tilde{L}\Sigma^{-1}L^{\top}\overline{U}_0^{\top}R\overline{U}_0).
		$
		For brevity, let $\overline{R} = \overline{U}_0^{\top}R\overline{U}_0$. Then, the Hessian satisfies
		\begin{align*}
		&\nabla^2 f(L,\Sigma)[(\tilde{L},\tilde{\Sigma}),(\tilde{L},\tilde{\Sigma})]\\
		& =\text{Tr}(\Sigma^{-1}\tilde{\Sigma}\Sigma^{-1}L^{\top}\overline{R}L\Sigma^{-1}\tilde{\Sigma}) +\text{Tr}(\Sigma^{-1}L^{\top}\overline{R}L \Sigma^{-1}\tilde{\Sigma}\Sigma^{-1} \tilde{\Sigma}) \\
		& ~~~+2\text{Tr}(\tilde{L}\Sigma^{-1}\tilde{L}^{\top}\overline{R})-4\text{Tr}(\overline{R}L\Sigma^{-1}\tilde{\Sigma}\Sigma^{-1}\tilde{L}^{\top})\\
		&=2\left\|R^{\frac{1}{2}}(\overline{U}_0\tilde{L}- \overline{U}_0L\Sigma^{-1} \tilde{\Sigma}) \Sigma^{-\frac{1}{2}}\right\|_F^2  \geq 0,
		\end{align*}
		which completes the proof.
	\end{proof}
	
	By Lemmas \ref{lem:asmp1} and \ref{lem:asmp2}, the LQR problem with covariance parameterization  \eqref{prob:equiV} and its convex reparameterization \eqref{prob:convex} satisfy the assumptions of \cite[Theorem 1]{sun2021analysis}. Then, there exists $c>0$ and a direction $Z\in \mathcal{N}(\overline{X}_0)$ with $\|Z\|_F=1$ in the descent cone of $\mathcal{S}(a)$ such that
	$
	J'(V)[Z] \leq -c (J(V)-J^*),
	$
	where $J'(V)[Z]$ denotes the derivative of $J(V)$ along the direction $Z$. Let 
	$Z' = \Pi_{\overline{X}_0} \nabla J(V) / \|\Pi_{\overline{X}_0} \nabla J(V)\|_F$ be the normalized projected gradient. Since both $Z$ and $Z'$ are feasible directions but $Z'$ is the direction of the projected gradient, we have $ J'(V)[Z'] \leq  J'(V)[Z], \forall Z\in \mathcal{N}(\overline{X}_0) $. Thus, it holds that $J(V) - J^* \leq \|\Pi_{\overline{X}_0} \nabla J(V)\|/c$.
	
	Next, we derive an explicit upper bound of $1/c$ over $V \in \mathcal{S}(a)$. By \cite[Theorem 1]{sun2021analysis}, $c$ is given by
	$$
	c = (2\max\{\|L-L^*\|_F/\underline{\sigma}(\Sigma), \|\Sigma - \Sigma^*\|_F  \|L\| /\underline{\sigma}^{2}(\Sigma) \})^{-1},
	$$
	where $(L^*,\Sigma^*)$ is an optimal solution of  (\ref{prob:convex}), and 
	$(L,\Sigma) = \arg\min_{(L',\Sigma') \in \mathcal{L}}\{ f(L',\Sigma'), \text{s.t.} ~L'(\Sigma')^{-1} = V\}.$
	We now provide upper bounds for $\underline{\sigma}^{-1}(\Sigma), \|\Sigma\|_F, \|L\|_F$. Since $\underline{\sigma}(\Sigma)\geq 1$, it holds $\underline{\sigma}^{-1}(\Sigma) \leq 1$. The sublevel set gives
	$\text{Tr}(Q\Sigma) \leq a$, and hence $\|\Sigma\|_F \leq a / \underline{\sigma}(Q)$. Since
	\begin{align*}
	&\underline{\sigma}(R) \underline{\sigma}^2(\overline{U}_0)\|\Sigma\|^{-1} \|L\|^2_F \leq \text{Tr}\left(L\Sigma^{-1}L^{\top}\overline{U}_0^{\top}R\overline{U}_0\right)  \\
	&\leq \text{Tr}(Q\Sigma) + \text{Tr}\left(L\Sigma^{-1}L^{\top}\overline{U}_0^{\top}R\overline{U}_0 \right)\leq a,
	\end{align*}
	an upper bound of $\|L\|_F$ is given by
	$$
	\|L\|_F \leq  \left(\frac{a\|\Sigma\|}{\underline{\sigma}(R) \underline{\sigma}^2(\overline{U}_0)}\right)^{{1}/{2}} \leq \frac{a}{(\underline{\sigma}(Q)\underline{\sigma}(R))^{1/2} \underline{\sigma}(\overline{U}_0)}.
	$$
	Those bounds are also true for $(L^*,\Sigma^*)$. Furthermore, we can provide the following upper bound 
	$$
	\frac{1}{c} \leq \frac{4a}{(\underline{\sigma}(Q)\underline{\sigma}(R))^{1/2} \underline{\sigma}(\overline{U}_0)} \max\left\{1, \frac{a}{\underline{\sigma}(Q)} \right\}.
	$$
	
	Since $\mathcal{S}(a)$ is non-empty, we have $a\geq J^* \geq \underline{\sigma}(Q)$. Hence,
	$$
	\frac{1}{c}\leq \frac{4a^2}{(\underline{\sigma}(Q))^{{3}/{2}}(\underline{\sigma}(R))^{{1}/{2}} \underline{\sigma}(\overline{U}_0)}:=\mu(a).
	$$

	\subsection{Proof of Lemma \ref{lem:smooth}}\label{app:2}
	We show the local smoothness of $J(V)$ by providing an upper bound for $\|\nabla^2 J(V)\| := \sup_{\|Z\|_F=1} \left|\nabla^2 J(V)[Z,Z]\right|$ over the sublevel set $V\in \mathcal{S}(a)$.
	We first present a formula for the Hessian of $J(V)$.
	\begin{lemma}\label{lem:hessian}
		For $V\in \mathcal{S}$ and a feasible direction $Z\in \mathbb{R}^{(n+m)\times n}$, the Hessian of $J(V)$ is characterized as
		\begin{align*}
		&\nabla^2 J(V)[Z,Z]=  4\text{Tr}\left(Z^{\top}\overline{X}_1^{\top}P'[Z]\overline{X}_1V\Sigma_V \right) \\
		& +2\text{Tr}\left(Z^{\top}\left(\overline{U}_0^{\top}R\overline{U}_0+\overline{X}_1^{\top}P_V\overline{X}_1\right)Z\Sigma_V\right),
		\end{align*}
		where
		$
		P'[Z] = \sum_{i=0}^{\infty} (V^{\top}\overline{X}_1^{\top})^i(Z^{\top}E_V + E_V^{\top}Z)(\overline{X}_1V)^i. 
		$
	\end{lemma}
	\begin{proof}
		By the definition of $P_V$ in (\ref{equ:Lya_P}), it holds that
		\begin{align*}
		P'[Z]& =Z^{\top}\overline{X}_1^{\top}P_V\overline{X}_1V+ V^{\top}\overline{X}_1^{\top}P_V\overline{X}_1Z \\
		& + Z^{\top}\overline{U}_0^{\top}R\overline{U}_0V + V^{\top}\overline{U}_0^{\top}R\overline{U}_0Z +V^{\top}\overline{X}_1^{\top}P'[Z]\overline{X}_1V  \\
		& = Z^{\top}E_V + E_V^{\top}Z + V^{\top}\overline{X}_1^{\top}P'[Z]\overline{X}_1V.
		\end{align*}
		
		Then, we have 
		$
		E'[Z] = (\overline{U}_0^{\top}R\overline{U}_0+\overline{X}_1^{\top}P_V\overline{X}_1)Z + \overline{X}_1^{\top}P'[Z]\overline{X}_1V,
		$
		and the Hessian $P''[Z]$ satisfies
		\begin{align*}
		P''[Z]&=  V^{\top}\overline{X}_1^{\top}P'[Z]\overline{X}_1Z +Z^{\top}\overline{X}_1^{\top}P'[Z]\overline{X}_1V  \\
		&~~~+Z^{\top}E'[Z] +(E'[Z])^{\top}Z +     V^{\top}\overline{X}_1^{\top}P''[Z]\overline{X}_1V\\
		&= S+V^{\top}\overline{X}_1^{\top}P''[Z]\overline{X}_1V,
		\end{align*}
		where 
		$
		S = 2V^{\top}\overline{X}_1^{\top}P'[Z]\overline{X}_1Z+2Z^{\top}\overline{X}_1^{\top}P'[Z]\overline{X}_1V + 2Z^{\top}(\overline{U}_0^{\top}R\overline{U}_0+\overline{X}_1^{\top}P_V\overline{X}_1)Z.
		$
		
		By the definition of $J(V)$ in (\ref{def:JG}), it holds that 
		$
		\nabla^2 J(V)[Z,Z]= \text{Tr} (P''[Z] )
		= \text{Tr}\left(S\Sigma_V\right).
		$ 
	\end{proof}
	
	Next, we provide an upper bound for $\|\nabla^2 J(V)\|$. Let $Z \in \mathbb{R}^{(n+m)\times n}$ be a feasible direction with $\|Z\|_F=1$. Then,
	\begin{equation}\label{equ:lemma7}
	\begin{aligned}
	&\|\nabla^2 J(V)\| \leq  4\left|\text{Tr}\left(Z^{\top}\overline{X}_1^{\top}P'[Z]\overline{X}_1V\Sigma_V \right)\right| \\ &+2\left\|Z^{\top}\left(\overline{U}_0^{\top}R\overline{U}_0+\overline{X}_1^{\top}P_V\overline{X}_1\right)Z\right\|\text{Tr}(\Sigma_V ).
	\end{aligned}
	\end{equation}

	For the first term, we have that
	\begin{align*}
	&\left|\text{Tr}(Z^{\top}\overline{X}_1^{\top}P'[Z]\overline{X}_1V\Sigma_V)\right|\\
	& \leq  \left\|Z^{\top}\overline{X}_1^{\top}P'[Z]\overline{X}_1V\Sigma_V^{\frac{1}{2}}\right\|_F\left\|\Sigma_V^{\frac{1}{2}}\right\|_F\\
	& \leq  \left\|\overline{X}_1\right\|_F^2  \left\|\overline{X}_1V\Sigma_V^{\frac{1}{2}}\right\|_F\left\|\Sigma_V^{\frac{1}{2}}\right\|_F \left\|P'[Z]\right\|_F \\
	&\leq \left\|\overline{X}_1\right\|_F^2  \|P'[Z]\|_F\cdot a/{\underline{\sigma}(Q)},
	\end{align*}
	where the last inequality follows from Lemma \ref{lem:bounds} and
	\begin{align*}
	&\left\|\overline{X}_1V\Sigma_V^{\frac{1}{2}}\right\|_F^2 
	=\text{Tr}\left(\overline{X}_1V\Sigma_VV^{\top}\overline{X}_1^{\top}\right)\\
	&\leq\text{Tr}\left(\overline{X}_1V\Sigma_VV^{\top}\overline{X}_1^{\top} + I_n\right) =  \left\|\Sigma_V^{\frac{1}{2}}\right\|_F^2.
	\end{align*}
	
	Thus, we only need to provide an upper bound for $\|P'[Z]\|_F$. Since $\|Z\|_F = 1$ and
	\begin{align*}
	Z^{\top}E_V + E_V^{\top}Z 
	&\preceq Z^{\top}\left(\overline{U}_0^{\top}R\overline{U}_0+\overline{X}_1^{\top}P_V\overline{X}_1\right)Z \\
	&~~~+ V^{\top}\left(\overline{U}_0^{\top}R\overline{U}_0+\overline{X}_1^{\top}P_V\overline{X}_1\right)V\\
	&=Z^{\top}\left(\overline{U}_0^{\top}R\overline{U}_0+\overline{X}_1^{\top}P_V\overline{X}_1\right)Z + P_V-Q \\
	&\preceq \left(\left(\xi(a)+a\right)/{\underline{\sigma}(Q)}-1\right)Q,
	\end{align*}
	the definition implies that 
	$
	P'[Z] = \sum_{i=0}^{\infty} (V^{\top}\overline{X}_1^{\top})^i(Z^{\top}E_V + E_V^{\top}Z)(\overline{X}_1V)^i 
	 \preceq ((\xi(a)+a)/{\underline{\sigma}(Q)}-1)P_V.
	$ 
	
	Then, it holds that $\|P'[Z]\|_F \leq a\left(\left(\xi(a)+a\right)/{\underline{\sigma}(Q)}-1\right)$, and the first term of the right-hand side of (\ref{equ:lemma7}) is bounded by 
	$4\left(\left(\xi(a)+a\right)/{\underline{\sigma}(Q)}-1\right) \|\overline{X}_1\|_F^2 {a^2}/{\underline{\sigma}(Q)}$. By Lemma \ref{lem:bounds}, the second term is bounded by
	$2 \xi(a)  {a}/{\underline{\sigma}(Q)}$. That is,
	\begin{align*}
	\|\nabla^2 J(V)\| \leq 4a^2\left(\xi(a) + a -\underline{\sigma}(Q)\right) \frac{\|\overline{X}_1\|_F^2}{\underline{\sigma}^2(Q)} +  \frac{2\xi(a) a}{\underline{\sigma}(Q)}.
	\end{align*}
	By the Taylor expansion of $J(V)$, the proof is completed.
	
	\subsection{Proof of Theorem \ref{thm:conv}}
	Let $V_\eta := V- \eta \Pi_{\overline{X}_0}\nabla J(V)$. We first show that for a non-optimal $V\in\mathcal{S}(a)$ and any $\eta \in [0,1/l(a)]$, it holds $V_{\eta} \in \mathcal{S}(a)$. 
	
	Define $\mathcal{S}^o(a):= \{V\in \mathcal{S} \mid J(V)<a\}$, and $(\mathcal{S}^{o}(a))^c$ as its complement, which is closed. By Lemma \ref{lem:smooth}, given $\phi\in (0,1)$, there exists $b>0$ such that $\|\nabla^2 J(V)\|\leq (1+\phi)l(a)$ for $V \in \mathcal{S}(a+b)$. Clearly, $\mathcal{S}(a) \cap (\mathcal{S}^{o}(a+b))^c = \emptyset$. Then, the distance between them $d:=\inf\{\|V'-V\|, \forall V \in \mathcal{S}(a), V'\in (\mathcal{S}^{o}(a+b))^c \}$ is positive. Let $\overline{N} \in \mathbb{N}_+$ be large enough such that $2/(\overline{N}(1+\phi)l(a)) < d/\|\Pi_{\overline{X}_0}\nabla J(V)\|$, which is well-defined since $V$ is not optimal. Define a stepsize $\tau \in [0, 2/(\overline{N}(1+\phi)l(a))]$. Since $\tau < d/\|\Pi_{\overline{X}_0}\nabla J(V)\|$, we have $\|V_{\tau}-V\|< d$, i.e., $V_{\tau}\in \mathcal{S}(a+b)$. Thus, we can apply Lemma \ref{lem:smooth} over $\mathcal{S}(a+b)$ to show
	$
	J(V_{\tau}) -  J(V) 
	\leq - \tau(1 -{(1+\phi)l(a)\tau}/{2})\| \Pi_{\overline{X}_{0}} \nabla J(V)\|^2\leq 0,
	$
	where the last inequality follows from $\tau \leq 2/(1+\phi)l(a)$. This implies that the segment between $V$ and $V_{\tau}$ is contained in $\mathcal{S}(a)$. It is also clear that $V_{2\tau}\in \mathcal{S}(a+b)$ since $\|V_{2\tau}-V_{\tau}\|< d$. Then, we can use induction to show that the segment between $V$ and $V_{N\tau}$ for $N\in \mathbb{N}_+$ is in $ \mathcal{S}(a)$ as long as $N\tau \leq 2/(1+\phi)l(a)$. Since $\phi \in (0,1)$, we let $\eta \leq 1/l(a)$ to ensure the segment between $V$ and $V_{\eta}$  contained in $ \mathcal{S}(a)$. 
	
	Then, a simple induction leads to that for $\eta \in [0,1/l^0]$, the update (\ref{equ:gd}) satisfies $V^k \in \mathcal{S}(J(V^0)), \forall k \in \mathbb{N}$. Moreover, 
	$$
	J(V^{k+1}) \leq  J(V^k) - \eta\left(1 -\frac{l^0\eta}{2}\right)\| \Pi_{\overline{X}_0} \nabla J(V^k)\|^2.
	$$	
	Using Lemma \ref{lem:pl} and subtracting $J^*$ in both sides yield
	$$
	J(V^{k+1})-J^* \leq J(V^k) - J^* - \frac{2\eta - l^0\eta^2}{2(\mu^0)^2}(J(V^k) - J^*)^2.
	$$ 
	Dividing by $(J(V^k) - J^*)(J(V^{k+1}) - J^*)$ in both sides and noting $J(V^{k+1})\leq J(V^{k}) $ lead to
	$$
	\frac{2\eta - l^0\eta^2}{2(\mu^0)^2} \leq \frac{1}{J(V^{k+1}) - J^*} - \frac{1}{J(V^k) - J^*}.
	$$
	Summing up both sides over $0,1,\dots,k-1$ and using telescopic cancellation yield that
	$$
	\frac{k(2\eta - l^0\eta^2)}{2(\mu^0)^2} \leq \frac{1}{J(V^k) - J^*}.
	$$ 
	Letting the right-hand side of the above inequality equal $\epsilon$ and solving $k$ yields the required number of iterations (\ref{equ:kk}).

	\section{Proof in Section \ref{sec:adap}}
	
	\subsection{Proof of Lemma \ref{lem:Ut}}
	
	Since
	$
	\underline{\sigma}(\overline{U}_{0,t}) \geq \underline{\sigma}(U_{0,t})\underline{\sigma}(D_{0,t})/t
	$, we provide lower bounds for $\underline{\sigma}(U_{0,t})$ and $\underline{\sigma}(D_{0,t})$, respectively.
	
	Let $\theta_l\in \mathbb{R}^m,\theta_r\in \mathbb{R}^t$ be the left and right singular vectors of $U_{0,t}$ corresponding to the minimal singular value $\underline{\sigma}(U_{0,t})$, respectively, i.e., $\theta_l^{\top}U_{0,t} = \underline{\sigma}(U_{0,t}) \theta_r^{\top}$ with $\|\theta_l\|=\|\theta_r\|=1$. Define the diagonal block matrix $\Omega:= \text{diag}(\theta_l^{\top}, \dots, \theta_l^{\top}) \in \mathbb{R}^{n\times nm}$, which has full row rank. Then, it holds that $\Omega\cdot \mathcal{H}_{n+1}(U_{0,t}) = \underline{\sigma}(U_{0,t}) \mathcal{H}_{n+1}(\theta_r)$, and hence
	\begin{equation}\label{equ:msv_U}
	\underline{\sigma}(\Omega) \underline{\sigma}(\mathcal{H}_{n+1}(U_{0,t})) \leq \underline{\sigma}(U_{0,t}) \|\mathcal{H}_{n+1}(\theta_r)\|.
	\end{equation}
	Since $\|\theta_l\|=\|\theta_r\|=1$, we have  $\|\mathcal{H}_{n+1}(\theta_r)\|\leq \sqrt{n+1}$ and $\underline{\sigma}(\Omega)=1$. Substituting them into \eqref{equ:msv_U} and using Assumption \ref{ass:PE}, we yield
	$
	\underline{\sigma}(U_{0,t}) \geq \sqrt{t}\gamma.
	$

	By Assumption \ref{ass:PE} and the robust fundamental lemma \cite[Theorem 5]{coulson2022quantitative}, there exist positive constants $\zeta, \omega$ such that
	$$
	\underline{\sigma}(D_{0,t}) \geq \sqrt{t}\gamma \zeta - \frac{\omega\|\mathcal{H}_n(W_{0,t-1})\|}{\sqrt{n+1}}.
	$$
	By Assumption \ref{ass:noise}, it holds that
	$
	\|\mathcal{H}_n(W_{0,t-1})\|\leq\|\mathcal{H}_n(W_{0,t-1})\|_F \leq \sqrt{n(t-n)}\delta\leq \sqrt{(n+1)t}\delta.
	$
	Then, the hypothesis ${\delta}/{\gamma} \leq \zeta/(2\omega)$ leads to $\underline{\sigma}(D_{0,t}) \geq \sqrt{t}\gamma \zeta/2$. 
	
	Hence, we have $\underline{\sigma}(\overline{U}_{0,t}) \geq \underline{\sigma}(U_{0,t})\underline{\sigma}(D_{0,t})/t \geq \gamma \zeta/2$.

	\subsection{Proof of Lemma \ref{lem:cost_diff}}
	We first provide the following result on the perturbation theory for Lyapunov equations.
	\begin{lemma}\label{lem:perturb}
		Let $A\in \mathbb{R}^{n \times n}$ be  stable  and $\Sigma(A) = I_n + A\Sigma(A) A^{\top}$. If 
		$
		\|A'-A\| \leq 1/({4\|\Sigma(A)\|(1+\|A\|)}),
		$
		then $A'$ is stable and
		$
		\|\Sigma(A')-\Sigma(A)\| \leq 4 \|\Sigma(A)\|^2(1+\|A\|)\|A'-A\| .
		$
	\end{lemma}
	\begin{proof}
		The proof follows the same vein as that of \cite[Lemma 16]{fazel2018global} and is omitted for saving space.
	\end{proof}
	
	1) Proof of \eqref{equ:1}. By definition, the cost can be written as
	\begin{align*}
	& J_{t}(V_t') = \text{Tr}\left((Q+K_{t}^{\top}RK_{t} )\Sigma_t(V_t')\right) \\
	& C(K_t) = \text{Tr}\left((Q+K_{t}^{\top}RK_{t} )\Sigma_{K_t}\right)
	\end{align*}
	with $\Sigma_t(V_t') = I_n + \overline{X}_{1,t}V_t'\Sigma_t(V_t')(V_t')^{\top} \overline{X}_{1,t}^{\top}$ and $\Sigma_{K_t} = I_n + (A+BK_t)\Sigma_{K_t} (A+BK_t)^{\top}$.
	Then, it holds that
	$
	|C(K_{t}) - J_t(V_t')| =|\text{Tr}(Q+K_t^{\top}RK_t)(\Sigma_{K_t} - \Sigma_t({V_t'}))|
	\leq J_t(V_t')\|\Sigma_{K_t} - \Sigma_t({V_t'})\|.$
	
	Next, we provide an upper bound for $\|\Sigma_{K_t} - \Sigma_t({V_t'})\|$ based on Lemma \ref{lem:perturb}.	Starting from (\ref{equ:diffff}), we have
	\begin{align*}
	&\left\|A+BK_t - \overline{X}_{1,t}V_t'\right\| 
	\leq \frac{2\delta}{\gamma \zeta}\left(1+ \sqrt{\frac{J_t(V_t')}{\underline{\sigma}(R)}} \right)\\ &\leq\frac{\underline{\sigma}^2(Q)}{4J_t(V_t')(\underline{\sigma}(Q)+J_t(V_t'))}
	\leq \frac{\underline{\sigma}(Q)}{4J_t(V_t')(1+\| \overline{X}_{1,t}V_t'\|)},
	\end{align*}
	where the first inequality follows from Lemmas \ref{lem:Ut} and \ref{lem:bounds}, the second from the condition on $\delta/\gamma$, and the last from  
	$
	\|\overline{X}_{1,t}V_t'\| \leq \|\overline{X}_{1,t}V_t'\Sigma_t^{{1}/{2}}(V_t')\|\|\Sigma_t^{{1}/{2}}(V_t')\| 
	\leq \|\Sigma_t(V_t')\| \leq {J_t(V_t')}/{\underline{\sigma}(Q)}.
	$ 
	
	Then, it follows from Lemma \ref{lem:perturb} that
	\begin{equation}\label{equ:sigma1}
	\begin{aligned}
	\left\|\Sigma_{K_t} - \Sigma_t(V_t')\right\| \leq& 4\left(\frac{J_{t}(V_t')}{\underline{\sigma}(Q)}\right)^2\left(1+\|\overline{X}_{1,t}V_t'\|\right)\\
	&\times \left\| A+BK_t  -\overline{X}_{1,t}V_t' \right\|.
	\end{aligned}
	\end{equation}
	
	Hence, the cost difference can be bounded as
	\begin{align*}
	&\left|C(K_{t}) - J_t(V_t')\right| 
	\leq J_t(V_t')\|\Sigma_{K_t} - \Sigma_t({V_t'})\| \\
	&\leq  \left(\frac{4J_{t}^3(V_t')}{\underline{\sigma}^2(Q)}\right)\left(1+\|\overline{X}_{1,t}V_t'\|\right)  \left\|A+BK_t - \overline{X}_{1,t}V_t'\right\| \\ 
	&\leq  {p(J_t(V_t'))\delta}/{(2\gamma)}.
	\end{align*}
	
	2) Proof of \eqref{equ:2}. By definition, we have
	$
	J_{t+1}(V_{t+1}) = \text{Tr}\left((Q+K_{t}^{\top}RK_{t} )\Sigma_{t+1}({V_{t+1}})\right) 
	$
	with  $\Sigma_{t+1}({V_{t+1}}) = I_n + \overline{X}_{1,t+1}V_{t+1}\Sigma_{t+1}({V_{t+1}})  V_{t+1}^{\top}\overline{X}_{1,t+1}^{\top}$. 
	Then, it follows that $
	|J_{t+1}(V_{t+1}) - J_{t}(V_t')| 
	= |\text{Tr}((Q+K_{t}^{\top}RK_{t})(\Sigma_{t+1}({V_{t+1}}) - \Sigma_t(V_t')))| 
	\leq J_t(V_t')\cdot \|\Sigma_{t+1}({V_{t+1}}) - \Sigma_t(V_t')\|.
	$
	Since
	\begin{equation}\label{equ:m1m23}
	\begin{aligned}
	&\left\| \overline{X}_{1,t+1}V_{t+1} -\overline{X}_{1,t}V_t' \right\|
	=\left\|\overline{W}_{0,t+1}V_{t+1} -  \overline{W}_{0,t}V_t'\right\|\\
	&= \left\| {W}_{0,t+1}D_{0,t+1}^{\dagger}\begin{bmatrix}
	K_{t} \\ I_n
	\end{bmatrix}  -  {W}_{0,t}D_{0,t}^{\dagger}\begin{bmatrix}
	K_{t} \\ I_n
	\end{bmatrix} \right\| \\
	& \leq \left\|{W}_{0,t+1}D_{0,t+1}^{\dagger} - {W}_{0,t}D_{0,t}^{\dagger}\right\| (1+\|K_{t}\|)\\
	& \leq \frac{4\delta}{\gamma \zeta}\left(1+ \sqrt{\frac{J_t(V_t')}{\underline{\sigma}(R)}} \right)\leq \frac{\underline{\sigma}(Q)}{4J_t(V_t')(1+\| \overline{X}_{1,t}V_t'\|)},
	\end{aligned}
	\end{equation}
	it follows from  Lemma \ref{lem:perturb} that
	\begin{align*}
	&\left\|\Sigma_{t+1}({V_{t+1}}) - \Sigma_t(V_t')\right\| \\
	&\leq 4\left(\frac{J_{t}(V_t')}{\underline{\sigma}(Q)}\right)^2\left(1+\|\overline{X}_{1,t}V_t'\|\right) \left\| \overline{X}_{1,t+1}V_{t+1} -\overline{X}_{1,t}V_t' \right\|\\
	&\leq \frac{16\delta}{\gamma \zeta}\left(\frac{J_{t}(V_t')}{\underline{\sigma}(Q)}\right)^2\left(1+\frac{J_{t}(V_t')}{\underline{\sigma}(Q)}\right)\left(1+ \sqrt{\frac{J_t(V_t')}{\underline{\sigma}(R)}} \right).
	\end{align*} 
	
	Then, the cost difference can be upper-bounded as
	$
	|J_{t+1}(V_{t+1}) - J_{t}(V_t')| \leq J_t(V_t') \|\Sigma_{t+1}({V_{t+1}}) - \Sigma_t(V_t')\| 
	 \leq p(J_t(V_t'))\delta/\gamma.
	$	
	
	3) Proof of \eqref{equ:3}. It follows by combining (\ref{equ:1}) and (\ref{equ:2}).
	
	\subsection{Proof of Lemma \ref{lem:opt_diff}}	
	Let $P^*$ and $P_t^*$ be the solution to the algebraic Riccati equation of the ground-truth $(A,B)$ and the estimated system $(\widehat{A}_t, \widehat{B}_t)$, respectively, i.e.,
	\begin{align*}
	&P^*=A^{\top}P^*A+Q- A^{\top} P^* B(R+B^{\top} P^* B)^{-1} B^{\top} P^* A,\\
	&P_t^*=\widehat{A}_t^{\top} P_t^* \widehat{A}_t+Q-\widehat{A}_t^{\top} P_t^* \widehat{B}_t(R+\widehat{B}_t^{\top} P_t^* \widehat{B}_t)^{-1} \widehat{B}_t^{\top} P_t^* \widehat{A}_t.
	\end{align*}
	
	By the perturbation theory for the discrete-time Riccati equation \cite[Proposition 1]{mania_certainty_2019}, there exist two constants $\bar{c}_1,\bar{c}_2$  depending only on $A,B,Q,R$, such that if $\|[\widehat{B}_t,\widehat{A}_t] - [B,A]\| \leq \bar{c}_1$, then
	$
	\|P_t^* - P^*\| \leq \bar{c}_2 \|[\widehat{B}_t,\widehat{A}_t] - [B,A]\|.
	$
	Since $\|[\widehat{B}_t,\widehat{A}_t]-[B,A]\| = \|W_{0,t}D_{0,t}^{\dagger}\| \leq 2\delta/(\gamma \zeta)$, it suffices to let $\delta/\gamma \leq  \zeta\bar{c}_1/2=:  \zeta c_1$. Then, it follows that
	$
	\|P_t^* - P^*\| \leq \bar{c}_2 \|[\widehat{B}_t,\widehat{A}_t] - [B,A]\|\leq 2\bar{c}_2 {\delta}/{(\gamma \zeta)}.
	$
	Hence, it holds that
	$
	\left|C_{\text{CE},t}^* - C^*\right| = |\text{Tr}(P_t^*)-\text{Tr}(P^*)| \leq 2\sqrt{n}\bar{c}_2\delta/(\gamma \zeta).
	$
	By Lemma \ref{lem:equi}, it holds that $|C^*-J_t^*| = |C^*-C_{\text{CE},t}^*|$. Letting $c_2 = 2\sqrt{n}\bar{c}_2$ completes the proof.
	
	\subsection{Proof of Lemma \ref{lem:boundcost}}
	The proof relies on the lower bounds of the gradient dominance and smoothness constants in Lemmas \ref{lem:uni_pl} and \ref{lem:smooth_V}. For brevity, let $\mu_{t}:= \mu(J_t(V_t))$ and  $\l_{t}: = l(J_t(V_t))$.
	\begin{lemma}\label{lem:lower_ul}
		For $V_t \in \mathcal{S}_t$, it holds that $l_{t}\geq c_3$ and $\mu_{t}\geq c_4$.
	\end{lemma}
	\begin{proof}
		By definition, it holds that $J_t(V_t) \geq J_t^* \geq \underline{\sigma}(Q)$. Then, $l_{t}$ and $\mu_{t}$ can be lower-bounded accordingly.
	\end{proof}
	
	We use mathematical induction to show $J_t(V_t)\leq \bar{J}$ for $t > t_0$. We first show that it holds at $t=t_0+1$. Since $J_{t_0}(V_{t_0}') = J_{t_0}^*\leq \bar{J}$, the assumption of Lemma \ref{lem:cost_diff} is satisfied, and 
	$
	J_{t_0+1}(V_{t_0+1}) \leq J_{t_0}(V_{t_0}') + {p(J_{t_0}(V_{t_0}') )\delta}/{\gamma \zeta}\leq  J_{t_0}^*+{\bar{p}\delta}/{\gamma}.
	$
	Using the condition $\delta/\gamma \leq {\eta }/{(2\bar{\mu}^2\bar{p})}$ and Lemma \ref{lem:lower_ul} lead to 
	$$
	J_{t_0+1}(V_{t_0+1}) \leq J_{t_0}^*+\frac{\eta}{2 \bar{\mu}^2} \leq  J_{t_0}^*+\frac{1}{2\bar{l} \bar{\mu}^2} \leq  J_{t_0}^*+\frac{1}{2c_3 c_4^2}\leq \bar{J}.
	$$
	
	Next, we assume that $J_t(V_t)\leq \bar{J}$ and show $J_{t+1}(V_{t+1})\leq \bar{J}$. Since $J_t(V_t)\leq \bar{J}$, we have $l_{t} \leq \bar{l}$ and  $\eta \leq 1/\bar{l} \leq 1/l_{t}$. Following the same vein as the proof of Theorem \ref{thm:conv}, we can show that $V_t'\in \mathcal{S}_{t}$ and 
	\begin{equation}\label{equ:sss}
	\begin{aligned}
	J_t(V_t') -J_t(V_t)&\leq   - \eta\left(1 -\frac{l_{t}\eta}{2}\right) \left\| \Pi_{\overline{X}_{0,t}} \nabla J_t(V_t)\right\|^2 \\
	&\leq  -\frac{\eta}{2} \left\| \Pi_{\overline{X}_{0,t}} \nabla J_t(V_t)\right\|^2,
	\end{aligned}
	\end{equation}
	where the second inequality follows from  $\eta \leq 1/l_t$. Hence, it holds that $J_t(V_t') \leq J_t(V_t) \leq \bar{J}$. By Lemma \ref{lem:cost_diff} and our condition $\delta/\gamma  \leq {\eta }/{(2\bar{\mu}^2\bar{p})}$, we have 
	\begin{equation}\label{equ:aaa}
	J_{t+1}(V_{t+1})  - J_t(V_t') \leq \frac{p(J_t(V_t'))\delta}{\gamma}\leq  \frac{\bar{p}\delta}{\gamma } \leq \frac{\eta}{2 \bar{\mu}^2}.
	\end{equation}
	
	Combing \eqref{equ:sss} and \eqref{equ:aaa} and using Lemma \ref{lem:uni_pl} yield
	\begin{equation}\label{equ:pp}
	\begin{aligned}
	J_{t+1}(V_{t+1}) - J_t(V_t)&\leq \frac{\eta}{2\bar{\mu}^2}- \frac{\eta}{2}\left\| \Pi_{\overline{X}_{0,t}} \nabla J_t(V_t)\right\|^2 \\
	& \leq  \frac{\eta}{2\mu_{t}^2} - \frac{\eta}{2\mu^2_{t}}(J_t(V_t) - J_t^*)^2.
	\end{aligned}
	\end{equation}
	
	To show $J_{t+1}(V_{t+1})\leq \bar{J}$, we consider two cases. Suppose that
	$J_t(V_t) \geq C^* + c_1c_2 + 1.$ By Lemma \ref{lem:opt_diff}, it holds that $ \left|C^* - J_t^*\right| \leq c_2\delta/(\gamma \zeta)
	\leq c_1c_2$. Hence, $
	(J_t(V_t) - J_t^*)^2 \geq (C^* + c_1c_2 + 1 - J_t^*)^2 \geq 1.
	$
	Furthermore, \eqref{equ:pp} yields
	$$
	J_{t+1}(V_{t+1}) \leq J_t(V_t)+  \frac{\eta}{2\mu_{t}^2} - \frac{\eta}{2\mu^2_{t}}(J_t(V_t) - J_t^*)^2  \leq \bar{J}.
	$$	
	Otherwise, if $J_t(V_t) < C^* + c_1c_2 +1,$ then \eqref{equ:pp} yields
	$$
	J_{t+1}(V_{t+1}) \leq J_t(V_t) 
	+ \frac{\eta}{2\mu_{t}^2}\leq J_t(V_t) + \frac{1}{2c_3c_4^2} \leq \bar{J},
	$$
	where the second inequality follows from $\eta \leq 1/\bar{l}$ and Lemma \ref{lem:lower_ul}. Thus, the induction shows that $J_t(V_t)\leq \bar{J}$ for $t > t_0$.
	
	It follows from \eqref{equ:sss} that $J_t(V_t')\leq \bar{J}$ for $t > t_0$. Recalling that $J_{t_0}(V_{t_0}') = J_{t_0}^*\leq \bar{J}$ completes the proof.
	
	\subsection{Proof of Theorem \ref{thm:final}}
	By Lemma \ref{lem:boundcost}, it holds that $V_t' \in \mathcal{S}_t, \forall t\geq t_0$. Hence, the assumption of Lemma \ref{lem:cost_diff} is satisfied and $\{K_t\}$ is stabilizing.
	
	It follows from $\eta \in (0, 1/\bar{l}]$ and Lemma \ref{lem:boundcost} that $\eta <1/l_{t+1}$. Then, the progress of the one-step projected gradient descent in \eqref{equ:M3} leads to
	$
	J_{t+1}(V_{t+1}') -J_{t+1}(V_{t+1}) \leq- ({\eta}/{2\mu_{t+1}^2})(J_{t+1}(V_{t+1}) - J_{t+1}^*)^2.
	$
	By Lemmas \ref{lem:cost_diff} and \ref{lem:boundcost}, it holds that
	$
	J_{t+1}(V_{t+1}) - J_t(V_t') \leq  {p(J_t(V_t'))\delta}/{\gamma} \leq  {\bar{p}\delta}/{\gamma}.
	$
	Adding the above two inequalities on both sides and noting $\mu_t\leq \bar{\mu}$ yield
	\begin{equation}\label{equ:kkk}
	J_{t+1}(V_{t+1}') - J_t(V_t') \leq  \frac{\bar{p}\delta}{\gamma} - \frac{\eta}{2\bar{\mu}^2}(J_{t+1}(V_{t+1}) - J_{t+1}^*)^2.
	\end{equation}
	
	By Lemma \ref{lem:cost_diff}, the left-hand side of \eqref{equ:kkk} is lower-bounded: 
	\begin{equation}\label{equ:ooo}
	J_{t+1}(V_{t+1}') - J_t(V_t') \geq C(K_{t+1}) - C(K_{t}) - \frac{\bar{p}\delta}{\gamma}.
	\end{equation}
	
	Next, we provide a lower bound for the right-hand side of \eqref{equ:kkk}. For $t > t_0$, it holds that
	\begin{align*}
	(J_t(V_t) - J_t^* )^2&= \left(J_t(V_t) - C^*\right)^2 + \left(C^* - J_t^* \right)^2\\
	& ~~~+ 2\left(J_t(V_t) - C^*\right)\left(C^* - J_t^* \right) \\
	&\geq \left(J_t(V_t) - C^*\right)^2 - 2\left|  J_t(V_t) - C^*\right|\left| C^* - J_t^*\right|\\
	&\geq \left(J_t(V_t) - C^*\right)^2 - 4c_2\bar{J}\delta/\gamma \zeta,
	\end{align*} 
	where the last inequality follows from Lemmas \ref{lem:opt_diff} and \ref{lem:boundcost}. 
	
	Define $e_t = C(K_t)-C^*$. Then, it follows that
	\begin{equation}\label{equ:43}
	\begin{aligned}
	&\left(J_t(V_t) - C^* \right)^2= \left(J_t(V_t) -  C(K_{t-1}) + e_{t-1} \right)^2 \\
	&\geq e_{t-1}^2 - 2\left| J_t(V_t) -  C(K_{t-1}) \right|\left|e_{t-1}\right|\\
	&\geq e_{t-1}^2 -\frac{3\bar{p}\delta}{\gamma} \left(C(K_{t-1}) + C^*\right)\\
	&\geq e_{t-1}^2 -\frac{3\bar{p}\delta}{\gamma } \left( J_t(V_t) + \frac{3\bar{p}\delta}{2\gamma} + C^*\right) \\
	&\geq e_{t-1}^2 -\frac{3\bar{p}\delta}{\gamma } \left( \bar{J} + \frac{3}{2}\left(C^* + \frac{1}{2 c_3c_4^2}\right) \right) \\
	&\geq e_{t-1}^2 - \frac{15\bar{p}\bar{J}\delta}{2\gamma},
	\end{aligned}
	\end{equation}
	where the second and the third inequalities follow from Lemma \ref{lem:cost_diff}, and the fourth follows from
	$$
	\bar{p}\frac{\delta}{\gamma} \leq \frac{\eta}{2 \bar{\mu}^2}\leq \frac{1}{2 \bar{l} \bar{\mu}^2}\leq \frac{1}{2 c_3c_4^2}.
	$$
	
	Then, the right-hand side of \eqref{equ:kkk} can be upper-bounded by
	$$
	- \frac{\eta}{2\bar{\mu}^2}e_{t}^2+ \left(\bar{p}+\frac{2\eta c_2\bar{J}}{\bar{\mu}^2\zeta} + \frac{15\eta \bar{p}\bar{J}}{4\bar{\mu}^2}\right) \frac{\delta}{\gamma }.
	$$
	
	Combining \eqref{equ:ooo}, we yield
	\begin{equation}\label{kkkk}
	e_{t+1} - e_{t} \leq - \frac{\eta}{2\bar{\mu}^2}e_{t}^2+ \left(2\bar{p}+\frac{2\eta c_2\bar{J}}{\bar{\mu}^2\zeta} + \frac{15\eta \bar{p}\bar{J}}{4\bar{\mu}^2}\right) \frac{\delta}{\gamma}.
	\end{equation}
	
	Let $\nu_3 =\sqrt{ 2\bar{\mu}^2( C(K_{t_0})-C^*)/\eta}$ and $\nu_4 = ({4\bar{p}\bar{\mu}^2}/{\eta} + 4c_2\bar{J}/\zeta+{15\bar{p}\bar{J}}/{2})^{1/2}$.
	Rearranging \eqref{kkkk} yields
	$
	e_{t}^2 \leq {2\bar{\mu}^2}(e_{t} - e_{t+1})/\eta + {\nu_4\delta}/{\gamma}.
	$
	Summing up both sides from $t_0$ to $t$ yields
	\begin{align*}
	\sum_{k=t_0}^{t}e_{k}^2 \leq \frac{2\bar{\mu}^2}{\eta}(e_{t_0} - e_{t+1}) + \frac{T\nu_4^2\delta}{\gamma }
	\leq \nu_3^2 + \frac{T\nu_4^2\delta}{\gamma}
	\end{align*}
	with $T=t-t_0+1$.	Dividing both sides by $T$ and using the AM-GM inequality yield that 
	$$
	\text{regret}_T \leq \sqrt{\frac{\nu_3^2}{T}+\frac{\nu_4^2\delta}{\gamma }}\leq \frac{{\nu_3}}{\sqrt{T}}  +  \nu_4\sqrt{\frac{ \delta}{\gamma }}.
	$$
	
	Letting $\nu_1 = 1/\bar{l}$ and $\nu_2$ be the right-hand side of \eqref{equ:cond_bound}, the proof is completed.
	
	
	\bibliographystyle{IEEEtran}
	\bibliography{mybibfile}

	\begin{IEEEbiography}[{\includegraphics[width=1in,height=1.25in,clip,keepaspectratio]{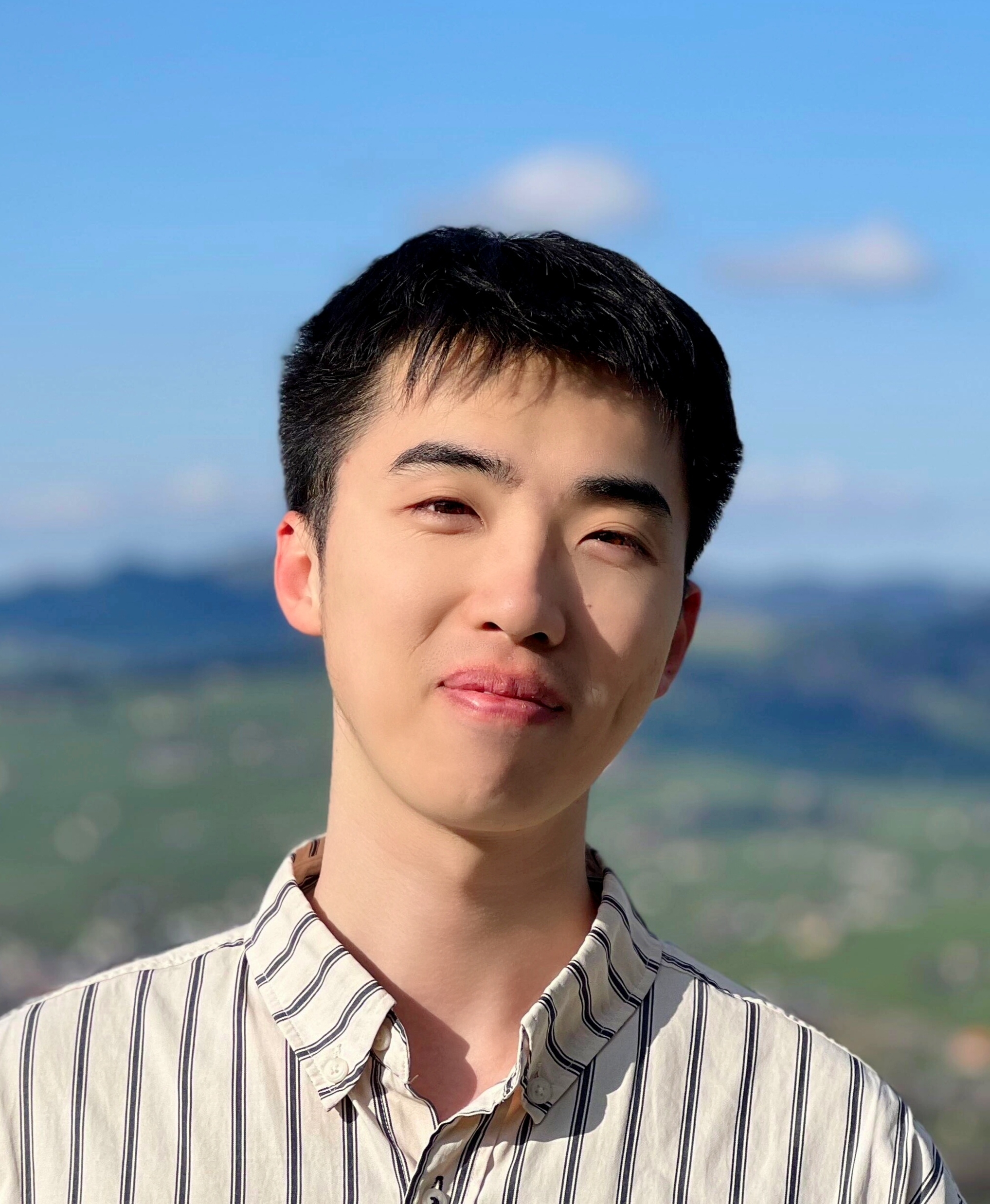}}]{Feiran Zhao} received the B.S. degree in  Control Science and Engineering from the Harbin Institute of Technology, China, in 2018, and the Ph.D. degree in Control Science and Engineering from the Tsinghua University, China, in 2024. He held a visiting position at ETH Z\"{u}rich. His research interests include policy optimization, data-driven control, adaptive control and their applications.
	\end{IEEEbiography}

	\begin{IEEEbiography}
		[{\includegraphics[width=1in,height=1.25in,clip,keepaspectratio]{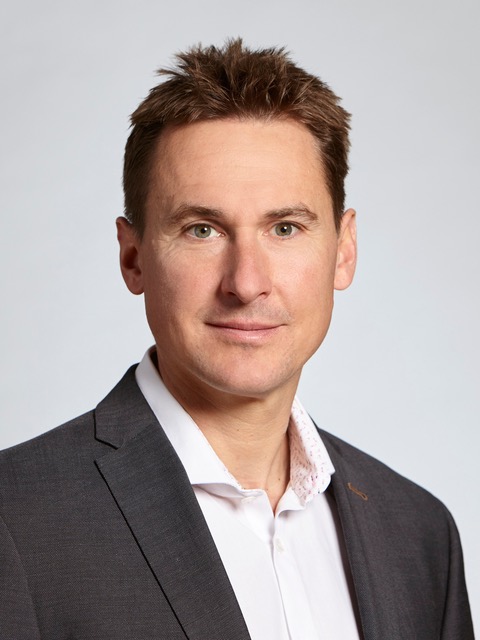}}]
%
		{Florian D\"{o}rfler} is a Full Professor at the Automatic Control Laboratory at ETH Z\"{u}rich. He received his Ph.D.
		degree in Mechanical Engineering from the University of California at Santa Barbara in 2013, and a
		Diplom degree in Engineering Cybernetics from the University of Stuttgart in 2008. From 2013 to 2014
		he was an Assistant Professor at the University of California Los Angeles. He has been serving as the
		Associate Head of the ETH Z\"{u}rich Department of Information Technology and Electrical Engineering 
		from 2021 until 2022. His research interests are centered around automatic control, system theory,
		and optimization. His particular foci are on network systems, data-driven settings, and applications to
		power systems. He is a recipient of the distinguished young research awards by IFAC (Manfred Thoma
		Medal 2020) and EUCA (European Control Award 2020). His students were winners or finalists for
		Best Student Paper awards at the European Control Conference (2013, 2019), the American Control
		Conference (2016, 2024), the Conference on Decision and Control (2020), the PES General Meeting
		(2020), the PES PowerTech Conference (2017), the International Conference on Intelligent Transportation Systems (2021), and the IEEE CSS Swiss Chapter Young Author Best Journal Paper Award (2022,
		2024). He is furthermore a recipient of the 2010 ACC Student Best Paper Award, the 2011 O. Hugo
		Schuck Best Paper Award, the 2012-2014 Automatica Best Paper Award, the 2016 IEEE Circuits and
		Systems Guillemin-Cauer Best Paper Award, the 2022 IEEE Transactions on Power Electronics Prize
		Paper Award, and the 2015 UCSB ME Best PhD award. He is currently serving on the council of the
		European Control Association and as a senior editor of Automatica.
	\end{IEEEbiography}

	\begin{IEEEbiography}
		[{\includegraphics[width=1in,height=1.25in,clip,keepaspectratio]{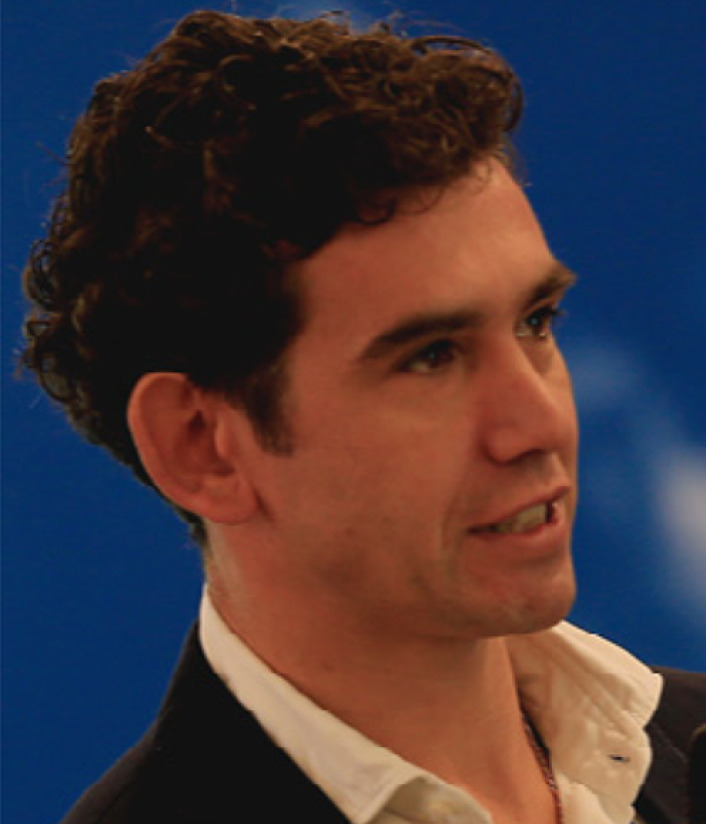}}]
		{Alessandro Chiuso} (Fellow, IEEE) is Professor with the Department of Information Engineering, Universit\`a di Padova. He received the ``Laurea" degree summa cum laude in Telecommunication Engineering from the University of Padova in July 1996 and the Ph.D. degree (Dottorato di ricerca) in System Engineering from the University of Bologna in 2000. He has been a visiting scholar with the Dept. of Electrical Engineering, Washington University St. Louis and Post-Doctoral fellow with the Dept. Mathematics, Royal Institute of Technology, Sweden. He joined the University of Padova as an Assistant Professor in 2001, Associate Professor in 2006 and then Full Professor since 2017. He currently serves as Editor (System Identification and Filtering) for Automatica. He has served as an Associate Editor for Automatica, IEEE Transactions on Automatic Control, IEEE Transactions on Control Systems Technology, European Journal of Control and MCSS. He is chair of the IFAC Coordinating Committee on Signals and Systems. He has been General Chair of the IFAC Symposium on System Identification, 2021 and he is a Fellow of IEEE (Class 2022). His research interests are mainly at the intersection of Machine Learning, Estimation, Identification and their Applications, Computer Vision and Networked Estimation and Control.
	\end{IEEEbiography}
	
	\begin{IEEEbiography}
		[{\includegraphics[width=1in,height=1.25in,clip,keepaspectratio]{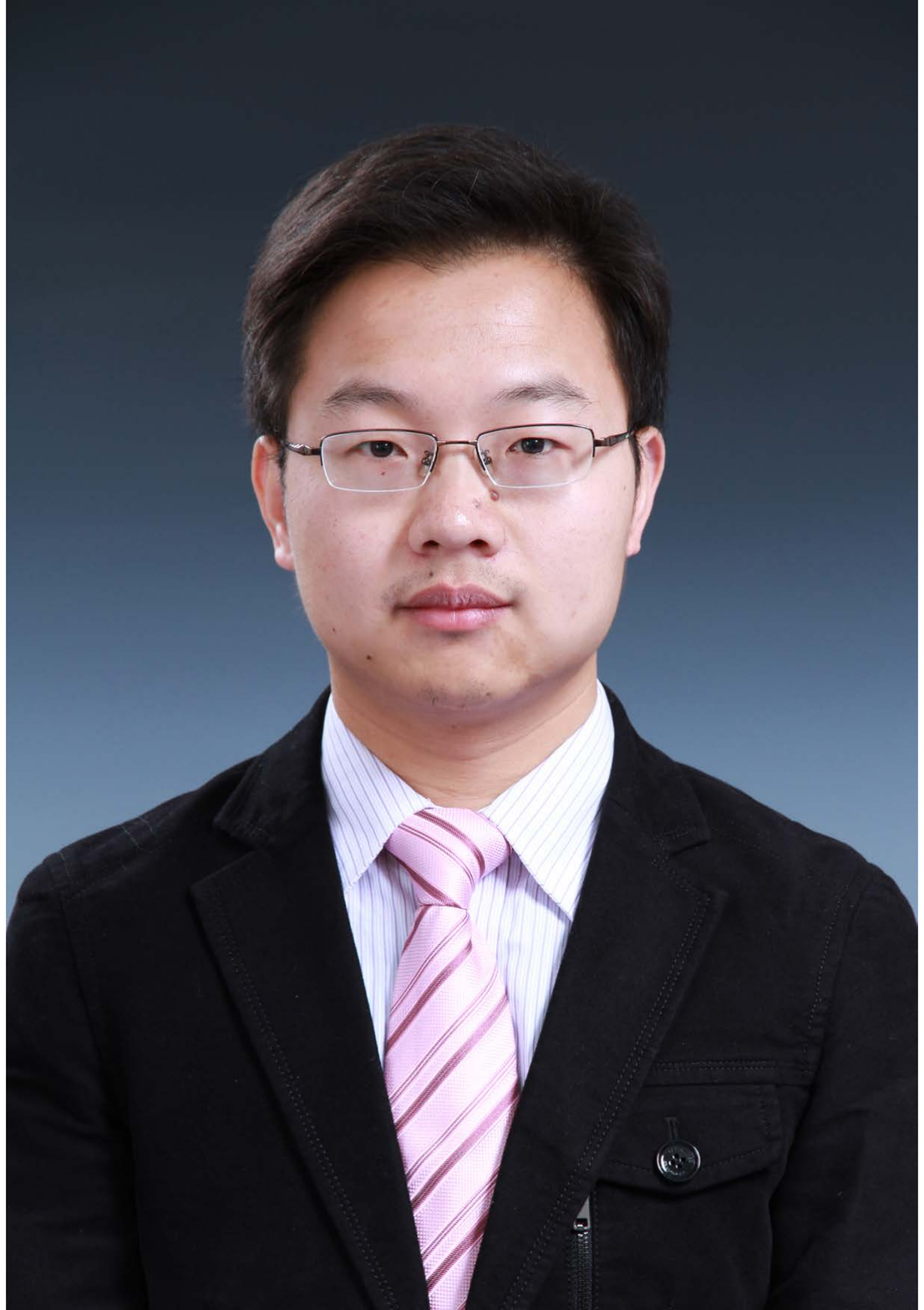}}]
		{Keyou You} received the B.S. degree in Statistical Science from Sun Yat-sen University, Guangzhou, China, in 2007 and the Ph.D. degree in Electrical and Electronic Engineering from Nanyang Technological University (NTU), Singapore, in 2012. After briefly working as a Research Fellow at NTU, he joined Tsinghua University in Beijing, China where he is now a Full Professor in the Department of Automation. He held visiting positions at Politecnico di Torino, Hong Kong University of Science and Technology, University of Melbourne and etc.
		
		Prof. You's research interests focus on the intersections between control, optimization and learning as well as their applications in autonomous systems. He received the Guan Zhaozhi award at the 29th Chinese Control Conference in 2010 and the ACA (Asian Control Association) Temasek Young Educator Award in 2019. He received the National Science Funds for Excellent Young Scholars in 2017, and for Distinguished Young Scholars in 2023. Currently, he is an Associate Editor for {\em Automatica}, {\em IEEE Transactions on Control of Network Systems}, and {\em IEEE Transactions on Cybernetics}.
	\end{IEEEbiography}

\end{document}